\documentclass[reqno,twoside]{amsart}
\usepackage{amscd}
\usepackage{amsfonts}
\usepackage{amsmath}
\usepackage{amssymb}
\usepackage{amsthm}
\usepackage{amsaddr}
\usepackage{fancyhdr}
\usepackage{latexsym}
\usepackage[colorlinks=true, pdfstartview=FitV, linkcolor=blue, citecolor=blue, urlcolor=blue]{hyperref}
\usepackage{enumitem}      
\usepackage{mathtools}            
\usepackage{indentfirst} 
\usepackage{color}
\usepackage{caption}          
\usepackage[normalem]{ulem}
\usepackage{upgreek}
\usepackage{dutchcal}
\usepackage{mathrsfs}
\newcommand{\eee}[1]{\begin{equation}#1\end{equation}}

\newcommand{\aaa}[1]{\begin{alignat}{2}#1\end{alignat}}

\newcommand{\nn}{\nonumber}
\newcommand{\p}{\partial}
\newcommand{\ve}{\varepsilon}

\newcommand{\no}[1]{\left\| #1 \right\|}

\usepackage{accents}

\usepackage{scalerel,stackengine}
\stackMath
\newcommand\wwhat[1]{%
\savestack{\tmpbox}{\stretchto{%
  \scaleto{%
    \scalerel*[\widthof{\ensuremath{#1}}]{\kern-.6pt\bigwedge\kern-.6pt}%
    {\rule[-\textheight/2]{1ex}{\textheight}}
  }{\textheight}%
}{0.5ex}}%
\stackon[1pt]{#1}{\tmpbox}%
}
\usepackage{thmtools}
\declaretheoremstyle[headfont=\normalfont\bfseries, bodyfont=\itshape, spaceabove=7pt, spacebelow=7pt]{theorem} 
\theoremstyle{theorem} 
\newtheorem{theorem}{Theorem}[] 
\newtheorem{remark}{Remark}[]

\newtheorem{corollary}{Corollary}[]

\theoremstyle{definition}

\allowdisplaybreaks
\numberwithin{equation}{section}
\usepackage{geometry}
\usepackage[breakable, theorems, skins]{tcolorbox}
\tcbset{enhanced}

\usepackage{tikz}
\usetikzlibrary{arrows.meta}
\usetikzlibrary{decorations.markings}
\tikzset{->-/.style={decoration={
  markings,
  mark=at position #1 with {\arrow{>}}},postaction={decorate}}}
  \tikzset{middlearrow/.style={
        decoration={markings,
            mark= at position 0.55 with {\arrow{#1}} ,
        },
        postaction={decorate}
    }
}
\makeatletter
  \renewcommand\subsection{\@startsection{subsection}{2}%
  \z@{-1\linespacing\@plus-0.7\linespacing}{0.7\linespacing}%
  {\bfseries}}
\makeatother
\let\OLDthebibliography\thebibliography
\renewcommand\thebibliography[1]{
  \OLDthebibliography{#1}
  \setlength{\parskip}{0pt}
  \setlength{\itemsep}{0pt plus 0.3ex}
}
\AtBeginDocument{
   \def\MR#1{}
}

\makeatletter
\@namedef{subjclassname@2020}{%
  \textup{2020} Mathematics Subject Classification}
\makeatother

\def\be{\begin{equation}}
	\def\ee{\end{equation}}
\def\bse{\begin{subequations}}
	\def\ese{\end{subequations}}

\def\sech{\mathop{\rm sech}\nolimits}

\def\N{\mathbb{N}}
\def\R{\mathbb{R}}

\def\eps{\epsilon}

\def\@#1{{\mathbf{#1}}}
\def\_#1{{\mathsf{#1}}}

\def\min{\mathop{\rm min}\nolimits}
\def\max{\mathop{\rm max}\nolimits}

\def\p{\partial}

\def\1{{\bf 1}}

\makeatother

\begin{document}

\title{On the proximal dynamics between integrable and non-integrable members of a generalized Korteweg-de Vries family of equations}

\author{Nikos I. Karachalios$^\dagger$, Dionyssios Mantzavinos$^*$, Jeffrey Oregero$^*$}
	
\address{
\normalfont $^\dagger$Department of Mathematics, University of Thessaly, Greece
\\
\normalfont $^*$Department of Mathematics, University of Kansas, USA
} 
\email{karan@uth.gr, mantzavinos@ku.edu, 
oregero@ku.edu \textnormal{(corresponding author)}}

\thanks{\textit{Acknowledgements.} DM gratefully acknowledges support from the U.S. National Science Foundation (NSF-DMS 2206270 and NSF-DMS 2509146) and the Simons Foundation (SFI-MPS-TSM-00013970). Furthermore, DM is thankful to the Department of Mathematics of the University of Thessaly, Lamia, Greece for their warm hospitality during February of 2025, when part of this work was undertaken.}
\subjclass[2020]{35Q53, 35B35, 37K40}
\keywords{Korteweg-de Vries (KdV) equation, modified KdV equation, generalized KdV equation, completely integrable systems, proximity estimates, soliton dynamics, structural stability} 
\date{November 18, 2025. \textit{Revised}: January 13, 2026}

\begin{abstract}
The distance between the solutions to the integrable Korteweg-de Vries (KdV) equation and a broad class of non-integrable generalized KdV (gKdV) equations is estimated in appropriate Sobolev spaces. This family of equations includes, as special cases, the standard gKdV equation with power nonlinearities as well as  weakly nonlinear perturbations of the KdV equation. For initial data and nonlinearity parameters of arbitrary size, we establish distance estimates based on a crucial size estimate for local gKdV solutions that grows linearly with the norm of the initial data. Consequently, these estimates predict that the dynamics of the gKdV and KdV equations remain close over long time intervals for initial amplitudes approaching unity, while providing an explicit rate of deviation for larger amplitudes. These theoretical results are supported by numerical simulations of one-soliton and two-soliton initial conditions, which show excellent agreement with the theoretical predictions. Furthermore, it is demonstrated that in the case of power nonlinearities and large solitonic initial data, the deviation between the integrable and non-integrable dynamics can be drastically reduced by incorporating suitable rotation effects via a rescaled KdV equation. As a result, the integrable dynamics stemming from the rescaled KdV equation may persist within the gKdV family of equations over remarkably long timescales.
\end{abstract}

\vspace*{-0.1cm}
\maketitle
\markboth
{Nikos I. Karachalios, Dionyssios Mantzavinos, Jeffrey Oregero}
{On the proximal dynamics between integrable and non-integrable members of a gKdV family of equations}


\section{Introduction and main results}
\label{s:intro}

We consider the family of generalized Korteweg-de Vries (gKdV) equations in one spatial dimension
\be
\label{gkdv}
U_t + U_{xxx} + F(U) \, U_x = 0, \quad x\in\mathbb R, \ t>0, 
\ee
where $U=U(x, t)$ is real-valued and $F:\R\to\R$ is a polynomial of degree $k$ satisfying $F(0)=0$, namely
\be
\label{e:polynon}
F(U) = \sum_{j=1}^k a_j U^j.
\ee
This class of nonlinearities includes the physically relevant case of the \textit{power nonlinearity},
\eee{\label{kdvk}
U_t + U_{xxx} + U^k U_x = 0, 
}
 which for $k=1$ corresponds to the celebrated KdV equation
\be
\label{kdv}
u_t + u_{xxx} + uu_{x} = 0
\ee
and for $k=2$ yields the modified KdV (mKdV) equation. The family \eqref{gkdv} also includes nonlinearities of the form $F(U) = U + \delta U^k$ with  $k\geq 2$ and $\delta\in\mathbb{R}$, which for $k=2$ and $\delta>0$ corresponds to the Gardner equation~\cite{gear1983second}. Importantly, this latter class of nonlinearities is strongly motivated by the seminal work \cite{dz2002}, which considers weak perturbations of the defocusing integrable cubic nonlinear Schr\"odinger (NLS) equation of the form
$iV_t + V_{xx} - 2|V|^2 V - \delta |V|^k V = 0$ for $k>2$ and small $\delta>0$. By analogy, the gKdV equation with $F(U) = U + \delta U^k$ and small $\delta$ corresponds to a weakly perturbed KdV equation.

The KdV, mKdV and Gardner equations  are prime examples of completely integrable systems. In particular, the KdV equation \eqref{kdv} possesses a rich collection of solutions, such as the famous one-soliton and multi-soliton solutions \cite{zk1965,akns1974,ablowitz1981solitons}, periodic time-translating solutions expressed in
terms of elliptic functions  known as cnoidal waves \cite{mckean1977stability} and their generalizations known as KdV wavetrains \cite{flaschka1980multiphase,mclaughlin1981modulations}, and even families of solutions that are periodic in the $x$-variable and almost periodic in $t$ lying in finite-dimensional tori  \cite{lax1975periodic,lax1976almost}.  The integrable members of the gKdV family have been studied via a plethora of techniques, including the powerful inverse scattering transform method \cite{ggkm1967,l1968,akns1974,miura1976korteweg,ablowitz1981solitons}. The Riemann-Hilbert problems arising through that process have been rigorously analyzed via the Deift-Zhou nonlinear steepest descent method, leading to explicit  asymptotics of the solution in the long-time limit \cite{dz1993}; we also refer to \cite{grunert2009long} for a survey of results regarding this method. With the exception of the KdV, mKdV and Gardner equations, the family of gKdV equations \eqref{gkdv} is not completely integrable. 

Regarding Hadamard well-posedness, namely the existence, uniqueness and continuous dependence of solution on the data, the Cauchy problem for the KdV equation on the infinite line with initial data in $L^2$-based Sobolev spaces $H^s(\R)$ has been studied extensively, see for example \cite{bs1975,st1976,kpv1989,kpv1991-kdv,b1993-kdv,kpv1996,ckstt2003} as well as the books \cite{t2006,lp2009} and the references therein. In fact, several of these works concern the case of power nonlinearity $U^k U_x$. 
In the more general case of the gKdV family \eqref{gkdv}, the well-posedness of the Cauchy problem on the infinite line was established in \cite{s1997}. Specifically, under appropriate conditions that contain the assumption~\eqref{e:polynon} considered in the present work, for initial data $U_0\in H^s(\R)$ with $s>\frac 12$ it was shown in Theorem~3.2 of \cite{s1997} that there exists $T=T(\no{U_0}_{H^s(\R)})>0$ (which grows to infinity as $U_0\to 0$) and a unique strong solution to the gKdV equation  \eqref{gkdv} such that $U\in C_t([0,T];H_x^s(\R))$.

The purpose of the present work is to rigorously study the dynamics generated by the integrable KdV equation against those of the gKdV family \eqref{gkdv}, which as noted earlier is in general not integrable. 
The study of how integrable structures persist under non-integrable perturbations is a central theme in the theory of Hamiltonian partial differential equations \cite{t2006}.  
As mentioned above, for the KdV equation, the integrable framework provides a profound understanding of soliton dynamics, multi-soliton interactions, and finite-gap quasi-periodic solutions.  
A natural and longstanding question is whether, and in what precise sense, these coherent structures and regular motions persist when integrability is broken, as in gKdV-type equations with higher-order nonlinearities. 
In this regard, we establish the following result on the distance of solutions between the integrable KdV equation and the gKdV family \eqref{gkdv}.
\begin{theorem}
\label{prox-t}
Given $k\in\N$ and $\epsilon>0$, consider the Cauchy problems for the gKdV and KdV equations \eqref{gkdv} and~\eqref{kdv} with initial data $U_0(x)$ and $u_0(x)$, respectively. For $s\in\N \cup \{0\}$, let  
\eee{\label{tcs}
T_s
= 
c_{s, k} 
\min\bigg\{
\Big(A_k \sum_{j=1}^k \no{U_0}_{H^{s+1}(\R)}^j \Big)^{-1},
\no{u_0}_{H^{s+1}(\R)}^{-1}
\bigg\}
}
where $c_{s, k} = c(s, k) > 0$ and $A_k := \max\left\{j |a_{j}|: 1\leq j \leq k\right\}$. 
\\[2mm]
\textnormal{(i) \underline{$L^2$ distance}:} If the initial data satisfy
\aaa{
\label{icorder}
&\no{u_0}_{H^2(\R)}\le c_0 \eps, \quad \no{U_0}_{H^2(\R)}\le C_0\eps, 
\\
&\no{u_0-U_0}_{L^2(\R)} \le C \max\left\{\eps^2, \eps^{k+1}\right\},
\label{icdorder}
}
for some constants $c_0$, $C_0$ and $C$, then for each fixed $T \in (0, T_1]$ there exists a constant $\widetilde{C} = \widetilde C(k, T)$ such that
\be\label{dl2}
\sup_{t\in[0,T]}\no{u(t)-U(t)}_{L_x^2(\R)} \le \widetilde{C}  \max\left\{\eps^2, \eps^{k+1}\right\}.
\ee
\textnormal{(ii) \underline{$H^s$ and $L^\infty$ distance for any $s\in \N$}:} If the initial data satisfy
\aaa{
\label{icorders}
&\no{u_0}_{H^{s+1}(\R)}\le c_{0, s} \eps, \quad \no{U_0}_{H^{s+1}(\R)}\le C_{0, s} \eps, 
\\
&\no{u_0-U_0}_{H^s(\R)} \le C_s \max\left\{\eps^2, \eps^{k+1}\right\},
\label{icdorders}
}
for some constants $c_{0, s}$, $C_{0, s}$ and $C_s$, then for each fixed $T \in (0, T_s]$ there exists a constant $\widetilde{C}_s = \widetilde C(s, k, T)$ such that
\be\label{dhs}
\sup_{t\in[0,T]}\no{u(t)-U(t)}_{L_x^\infty(\R)}
\leq
\sup_{t\in[0,T]}\no{u(t)-U(t)}_{H_x^s(\R)} \le \widetilde{C}_s  \max\left\{\eps^2, \eps^{k+1}\right\}.
\ee
\end{theorem}

In the special case of $F(U) = U^k$, for $s \in \N$, $s\neq 1$ and $U_0 \in H^s(\R)$, let $U\in C_t([0,T_0];H_x^s(\R))$ with $T_0=T_0(\no{U_0}_{H^s(\R)})>0$ be the unique strong solution  to the Cauchy problem for the gKdV equation with power nonlinearity \eqref{kdvk}. Then, Theorem \ref{prox-t} readily implies the following result.
\begin{corollary}[Power nonlinearity proximity]
\label{prox-c}
Given $k\in\N$ and $\epsilon>0$, consider the Cauchy problems for the KdV equation~\eqref{kdv} and the gKdV equation with a power nonlinearity \eqref{kdvk}, with initial data $u_0(x)$ and $U_0(x)$ respectively. For $s\in\N \cup \{0\}$, let  
\eee{\label{tcs-c}
\widetilde T_s
= 
c_{s, k} 
\min\left\{
\no{U_0}_{H^{s+1}(\R)}^{-k},
\no{u_0}_{H^{s+1}(\R)}^{-1}
\right\},
\quad 
c_{s, k} = c(s, k) > 0.
}
\noindent
\textnormal{(i) \underline{$L^2$ distance}:} If the initial data satisfy
\aaa{
\label{icorder-c}
&\no{u_0}_{H^2(\R)}\le c_0 \eps, \quad \no{U_0}_{H_x^2(\R)}\le C_0\eps, 
\\
&\no{u_0-U_0}_{L^2(\R)} \le C \max\left\{\eps^2, \eps^{k+1}\right\},
\label{icdorder-c}
}
for some constants $c_0$, $C_0$ and $C$, then for each fixed $T \in (0, \widetilde T_1]$ there exists a constant $\widetilde{C} = \widetilde C(k, T)$ such that
\be\label{dl3}
\sup_{t\in[0,T]}\no{u(t)-U(t)}_{L_x^2(\R)} \le \widetilde{C}  \max\left\{\eps^2, \eps^{k+1}\right\}.
\ee
\textnormal{(ii) \underline{$H^s$ and $L^\infty$ distance for any $s\in \N$}:} If the initial data satisfy
\aaa{
\label{icorders-c}
&\no{u_0}_{H^{s+1}(\R)}\le c_{0, s} \eps, \quad \no{U_0}_{H^{s+1}(\R)}\le C_{0, s} \eps, 
\\
&\no{u_0-U_0}_{H^s(\R)} \le C_s \max\left\{\eps^2, \eps^{k+1}\right\},
\label{icdorders-c}
}
for some constants $c_{0, s}$, $C_{0, s}$ and $C_s$, then for each fixed $T \in (0, \widetilde T_s]$ there exists a constant $\widetilde{C}_s = \widetilde C(s, k, T)$ such that
\be\label{dhs-c}
\sup_{t\in[0,T]}\no{u(t)-U(t)}_{H_x^s(\R)} \le \widetilde{C}_s  \max\left\{\eps^2, \eps^{k+1}\right\}.
\ee
Consequently, by the Sobolev embedding theorem, there exists constant $\widetilde C_s'=\widetilde C'(s, k, T)$ such that
\be\label{dsup-c}
\sup_{t\in[0,T]}\no{u(t)-U(t)}_{L_x^\infty(\R)} \le \widetilde{C}_s'  \max\left\{\eps^2, \eps^{k+1}\right\}.
\ee
\end{corollary}

\begin{remark}
Via the triangle inequality, the results of Theorem \ref{prox-t} and Corollary \ref{prox-c} readily imply distance estimates between \textit{any two members} of the gKdV family \eqref{gkdv}.
\end{remark}

Key to the proof of Theorem~\ref{prox-t} is the following size estimate  for solutions to the gKdV family \eqref{gkdv}.

\begin{theorem}
\label{se-t}
Suppose $s \in \N$ with $s\neq 1$ and let $F$ satisfy \eqref{e:polynon}. For  initial data $U_0 \in H^s(\R)$, let $U\in C_t([0,T_0];H_x^s(\R))$ be the unique strong solution  to the Cauchy problem of the gKdV equation \eqref{gkdv}, as guaranteed by  Theorem 3.2 in \cite{s1997}, where $T_0=T_0(\no{U_0}_{H^s(\R)})>0$. 
 Then, for 
\eee{\label{lspan}
T 
:=
\frac{c_{s,k}}{\no{U_0}_{H^s(\R)}^k + A_k \sum_{j=1}^{k-1} \no{U_0}_{H^s(\R)}^j} \leq T_0
}
with $A_k := \max\left\{j |a_{j}|: 1\leq j \leq k-1\right\}$, the solution  admits the estimate
\be
\label{se}
\no{U(t)}_{H_x^s(\R)} \le 2\no{U_0}_{H^s(\R)}, \quad 0 \leq t \leq T.
\ee
\end{theorem}

The bound \eqref{se} provides the analogue, in the case of the gKdV family \eqref{gkdv}, of the a priori bound given on page~51 of~\cite{et2016book} for the KdV equation \eqref{kdv}. It allows us to prove Theorem~\ref{prox-t} by estimating the equation satisfied by the difference of solutions $u-U$ to the KdV and gKdV equations and, more specifically, to derive bounds for $u-U$ not only in the $H^1(\mathbb{R})$ norm (as for the NLS-type equations studied in \cite{hkmms2024}) but also in the higher-regularity setting of $H^s(\mathbb{R})$ with $s\in\mathbb N$.

In the special case of gKdV with a power nonlinearity \eqref{kdvk},  $A_k = 0$ since $a_j=0$ for all $0\leq j\leq k-1$, thus  Theorem~\ref{se-t} readily yields the following result.
\begin{corollary}[Power nonlinearity size estimate]
\label{c:apriori}
Suppose $s \in \N$, $s\neq 1$, and $F(U) = U^k$. For initial data $U_0 \in H^s(\R)$, let $U\in C_t([0,T_0];H_x^s(\R))$ with $T_0=T_0(\no{U_0}_{H^s(\R)})>0$ be the unique strong solution  to the Cauchy problem of the gKdV equation with power nonlinearity \eqref{kdvk} guaranteed by  Theorem 3.2 of \cite{s1997}. 
 Then, for 
\eee{\label{lspank}
T 
=
\frac{c_{s,k}}{\no{U_0}_{H^s(\R)}^k}>0
}
the solution admits the  estimate
\be
\label{sec}
\no{U(t)}_{H_x^s(\R)} \le 2\no{U_0}_{H^s(\R)}, \quad 0 \leq t \leq T.
\ee
\end{corollary}

The interpretation of Theorem~\ref{prox-t} in terms of the persistence of integrable KdV dynamics in gKdV equations is illustrated by examining soliton propagation and interactions, particularly in the small data regime $0<\varepsilon<1$ where proximal dynamics between the two systems are expected.
When the KdV ``reference solution'' $u$ represents a one-soliton or a multi-soliton, Theorem~\ref{prox-t} implies that the gKdV solution $U$ remains $\mathcal{O}(\varepsilon^2)$ close in $H^s(\mathbb{R})$ for $t\in[0,T_s]$, guaranteeing structural persistence over a potentially long, but finite, time interval. Consequently, the gKdV soliton profiles, amplitudes and phases experience only slow deviations of order $\mathcal{O}(\varepsilon^2)$ from the corresponding ones for KdV, primarily due to minor radiative losses. Although non-integrable dynamics may eventually deviate from the exact KdV soliton solutions, Theorem~\ref{prox-t}  establishes a \textit{metastable regime} in which the interaction structure and solitary wave character are accurately retained. This fact is reinforced by the linear in time growth estimate (see inequality \eqref{lgc} in the proof of Theorem~\ref{prox-t})
\be\label{lgr}
\sup_{t\in[0,T]}\no{u(t)-U(t)}_{H_x^s(\R)} \lesssim \mathcal{O}\big(\max\{\varepsilon^2, \varepsilon^{k+1}\} t\big).
\ee

Our numerical simulations, which focus on the dynamics of both solitary wave and multi-soliton initial data defined by the corresponding analytical KdV solutions, confirm the above interpretation in excellent agreement with the theoretical results.
First, the numerically observed evolution of the distance $u(t)-U(t)$ (measured in $H_x^s(\mathbb{R})$ with $s=0,1,2$ and in  $L^\infty_x(\mathbb{R})$) exhibits linear growth in time, in remarkable agreement with the theoretical growth estimate of Theorem~\ref{prox-t}. Consequently, for small initial data and/or weak nonlinearities, the simulations confirm a remarkable persistence of the KdV soliton structure and its interactions within the gKdV flow. 

We find that the soliton profiles retain their coherence and velocity not only up to visible times far beyond the theoretical scale for small initial data but also for larger soliton amplitudes $\gtrsim 1$.
In paticular, for larger amplitudes, in the case of the power nonlinearity $F(U) = U^k$ we analyze in detail the phenomenon of prolonged proximal soliton dynamics relative to a \textit{rescaled} KdV equation, obtained via the transformation $t \mapsto t/\nu$ applied to the standard KdV equation~\eqref{kdv}. Assisted by numerical simulations of the gKdV soliton dynamics, we show that, for a given nonlinearity exponent $k \geq 2$, there exists a scaling parameter $\nu(k) < 1$ such that the gKdV soliton dynamics remains remarkably proximal to the analytical rescaled KdV soliton over significant time intervals. This effect, which geometrically corresponds to the \textit{rotation} of the linear caustic of the integrable KdV soliton, is analyzed  in detail. Quantitatively, the timescales required for the distance $u-U$ to reach $\mathcal{O}(1)$ are $\mathcal{O}(10^3)$ for $k=2$ (integrable mKdV) and $k=3$, decreasing to $\mathcal{O}(10^2)$ in the critical case $k=4$. Within these large timescales, the observed deviation is primarily due to weak radiative effects exhibited by the gKdV dynamics, while the soliton’s initial amplitude and shape persist with only minor differences from the analytical rescaled soliton profile.

At this point, it is important to highlight a fundamental distinction between the present approach and the one of \cite{dz2002}, particularly regarding weak nonlinearities. The methodology in \cite{dz2002} utilizes integrability theory by combining the inverse scattering transform for the analysis of the reflection coefficient in the defocusing NLS case with detailed estimates for the perturbed non-integrable defocusing NLS model. On the other hand, our approach for deriving distance estimates between the integrable and non-integrable models remains generic across the class of nonlinearities described by \eqref{e:polynon} for the (non-integrable) gKdV family of equations. Since our framework is independent of the size of the initial data and other parameters, it captures both proximity and deviation between the solutions of the compared models without explicitly relying on the integrability theory for the KdV equation. A potential synthesis of the methodology in \cite{dz2002} with the methods developed herein could yield even more precise distance estimates, presenting a formidable and challenging task for future investigations.

In summary, the rigorous analysis and numerical results present a coherent picture of long-time soliton persistence and near-integrable behavior in the gKdV dynamics. This work is organized as follows. In Section~\ref{s:apriori}, we prove Theorem \ref{se-t} for the size estimate of the gKdV solutions. In Section~\ref{s:proximity}, we establish the estimates of Theorem~\ref{prox-t} on the distance between the gKdV and KdV solutions. Finally, in Section~\ref{s:numerical}, we provide the results of our numerical simulations and conclude with a final discussion of the results, highlighting further potential applications.

\section{Proof of the size estimate}
\label{s:apriori}
In this section, we establish Theorem \ref{se-t} on the size estimate for the unique strong solution to the gKdV family~\eqref{gkdv} guaranteed by Theorem 3.2 of \cite{s1997} in the case of initial data in Sobolev spaces $H^s(\R)$. This size estimate will be crucial to the proximity analysis of Section \ref{s:proximity} that leads to the proof of Theorem \ref{prox-t}. In the case of smooth solutions, our result extends the estimate (3.4) derived in \cite{et2016book} for the KdV equation  to the gKdV setting. More generally, it extends that estimate to solutions in the Sobolev $H^s(\R)$ class. 
We first establish Theorem \ref{se-t} for data in the Schwartz space $\mathcal{S}(\R)$ of smooth and rapidly decaying functions. Then, we use a density argument to treat the case of Sobolev solutions.

In the case of $U_0 \in \mathcal S(\R)$, the results of \cite{k1983} (see also \cite{t1969,sjo1970,bs1975}) imply the existence of a unique solution in $U \in C_t^\infty([0, T_0]; \mathcal S(\R))$, where $T_0>0$ depends only on the norm of $U_0$. 
For that smooth solution, let $n\in\N\cup\{0\}$ with $0\le n\le s$ and apply $\p_x^n$ to equation \eqref{gkdv}.
Then, multiply the resulting equation by $\p_x^nU$ and integrate in $x$ while recalling that $U$ is real-valued to obtain
\be
\label{e:ident2}
\frac{1}{2} \frac{d}{dt} \no{\p_x^nU}_{L_x^2(\R)}^2 = 
-\int_{\R}\p_x^nU\p_x^n(F(U) \, U_x) \, dx,
\ee 
where $\int_{\R}\p_x^nU\p_x^{n+3}U dx=0$ via integration by parts.  

By the Leibniz rule for the $n$th derivative of a product,
\begin{align}
\int_{\R}\p_x^nU\p_x^n(F(U) \, \p_xU) dx 
&= \sum_{j=0}^n\tbinom{n}{j} \int_{\R} \p_x^nU\p_x^{n-j}(F(U)) \, \p_x^{j+1}U dx
\nn\\
&= \sum_{j=0}^{n-1}  \tbinom{n}{j} \int_{\R} \p_x^nU\p_x^{n-j}(F(U)) \, \p_x^{j+1}U dx + \int_{\R} (\p_x^nU)F(U) \, \p_x^{n+1}U dx 
\nn\\
&= \sum_{j=0}^{n-1}\tbinom{n}{j} \int_{\R} \p_x^nU\p_x^{n-j}(F(U)) \, \p_x^{j+1}U dx + \frac{1}{2}\int_{\R} F(U) \, \p_x(\p_x^nU)^2 dx. 
\end{align}
Thus, integrating by parts in the second integral and using the triangle inequality, we have
\eee{\label{injn}
\frac{1}{2} \frac{d}{dt} \no{\p_x^nU}_{L_x^2(\R)}^2
\leq |I_n| + \frac 12 |J_n|
}
where 
\eee{
I_n :=\sum_{j=0}^{n-1}\tbinom{n}{j} \int_{\R} \p_x^nU\p_x^{n-j}(F(U)) \, \p_x^{j+1}U dx,
\quad
J_n := \int_{\R} \p_x(F(U))\, (\p_x^nU)^2 dx.
}

Concerning $J_n$, we have
\aaa{
|J_n| 
\leq 
\int_{\R} \left|\p_x(F(U)) \, (\p_x^nU)^2\right| dx
&= \int_{\R} \left|F'(U) U_x (\p_x^nU)^2\right| dx
\nn\\
&\le  \no{F'(U)}_{L_x^\infty(\R)} \no{U_x}_{L_x^\infty(\R)} \no{\p^n_xU}_{L_x^2(\R)}^2.
}
Since $F'(U) = \sum_{j=1}^{k} j a_{j} U^{j-1}$, by the triangle inequality and the Sobolev embedding theorem,
\eee{
\no{F'(U)}_{L_x^\infty(\R)} 
\leq 
\sum_{j=1}^{k} j \left|a_{j}\right| \no{U^{j-1}}_{L_x^\infty(\R)}
\leq
c \sum_{j=1}^{k} j \left|a_{j}\right| \no{U^{j-1}}_{H_x^1(\R)}
}
so by the algebra property in $H^1(\R)$ we obtain
\eee{\label{f'b}
\no{F'(U)}_{L_x^\infty(\R)}
\leq
c \sum_{j=1}^{k} j \left|a_{j}\right| \no{U}_{H_x^1(\R)}^{j-1}
\leq
c \sum_{j=1}^{k} j \left|a_{j}\right| \no{U}_{H_x^s(\R)}^{j-1}
\leq
c \, \Big( k \left|a_k\right| \no{U}_{H_x^s(\R)}^{k-1} + A_k \sum_{j=1}^{k-1} \no{U}_{H_x^s(\R)}^{j-1} \Big)
}
where 
\eee{\label{akdef}
A_k :=  \underset{1\leq j \leq k-1}{\max} \left(j |a_{j}|\right).
} 
In turn, we find
\aaa{\label{jn}
|J_n| 
&\leq
c \, \Big( k \left|a_k\right| \no{U}_{H_x^s(\R)}^{k-1} + A_k \sum_{j=1}^{k-1} \no{U}_{H_x^s(\R)}^{j-1} \Big)
 \no{U_x}_{L_x^{\infty}(\R)}\no{U}_{H_x^n(\R)}^2
\nn\\
&\leq 
c \, \Big( k \left|a_k\right| \no{U}_{H_x^s(\R)}^{k-1} + A_k \sum_{j=1}^{k-1} \no{U}_{H_x^s(\R)}^{j-1} \Big)
\no{U}_{H_x^2(\R)}\no{U}_{H_x^n(\R)}^2
\nn\\
&\leq
c_k \Big(\no{U}_{H_x^s(\R)}^{k+2} + A_k \sum_{j=1}^{k-1} \no{U}_{H_x^s(\R)}^{j+2} \Big),
\quad 
c_k := c\max\{1, k \left|a_k\right|\},
}
with the last inequality in light of the fact that $s \geq \max\{2, n\}$. 

Regarding $I_{n}$, we first apply the Cauchy-Schwarz inequality to obtain
\eee{
|I_{n}| 
\leq 
\no{U}_{H_x^s(\R)}
\sum_{j=0}^{n-1}\tbinom{n}{j} 
\no{\p_x^{n-j}(F(U)) \, \p_x^{j+1}U}_{L_x^2(\R)}.
}
For $j=n-1$, as in the estimation of $J_n$, the bound \eqref{f'b} and the Sobolev embedding theorem yield
\aaa{\label{in1}
\no{\p_x(F(U)) \, \p_x^{n}U}_{L_x^2(\R)}
&= \no{F'(U) U_x \p_x^{n}U}_{L_x^2(\R)}
\le \no{F'(U)}_{L_x^\infty(\R)} \no{U_x}_{L_x^\infty(\R)} \no{\p_x^n U}_{L_x^2(\R)}
\nn\\
&\leq 
c\, \Big( k \left|a_k\right| \no{U}_{H_x^s(\R)}^{k+1} + A_k \sum_{j=1}^{k-1} \no{U}_{H_x^s(\R)}^{j+1} \Big)
\nn\\
&\leq 
c_k \Big(\no{U}_{H_x^s(\R)}^{k+1} + A_k \sum_{j=1}^{k-1} \no{U}_{H_x^s(\R)}^{j+1} \Big).
}
For $0\le j\le n-2$, by the Sobolev embedding theorem,
\eee{
\no{\p_x^{n-j}(F(U)) \, \p_x^{j+1}U}_{L_x^2(\R)}
\le \no{\p_x^{j+1}U}_{L_x^{\infty}(\R)} \no{\p_x^{n-j}(F(U))}_{L_x^2(\R)}
\le c \no{U}_{H_x^{j+2}(\R)} \no{\p_x^{n-j}(F(U))}_{L_x^2(\R)}
}
so, by the fact that $j+2 \leq n \leq s$, 
\eee{
\no{\p_x^{n-j}(F(U)) \, \p_x^{j+1}U}_{L_x^2(\R)}
\le c \no{U}_{H_x^s(\R)} \no{\p_x^{n-j}(F(U))}_{L_x^2(\R)}.
}
Moreover, by the triangle inequality and the algebra property, 
\aaa{
\no{\p_x^{n-j} (F(U))}_{L_x^2(\R)}
&\leq
\no{F(U)}_{H_x^{n-j}(\R)}
\leq
\Big\|\sum_{j=1}^k a_j U^j\Big\|_{H_x^s(\R)}
\nn\\
&\leq
\sum_{j=1}^k |a_j| \no{U^j}_{H_x^s(\R)}
\leq
c \sum_{j=1}^k |a_j| \no{U}_{H_x^s(\R)}^j
\nn\\
&\leq
c\, \Big( |a_k| \no{U}_{H_x^s(\R)}^k + B_k \sum_{j=1}^{k-1} \no{U}_{H_x^s(\R)}^{j} \Big)
\label{e:chainrule_est}
}
where  \eee{\label{bkdef}B_k := \underset{1\leq j \leq k-1}{\max} |a_j|.}
Therefore,  
\eee{\label{in2}
\no{\p_x^{n-j}(F(U)) \, \p_x^{j+1}U}_{L_x^2(\R)}
\le 
c\, \Big( |a_k| \no{U}_{H_x^s(\R)}^{k+1} + B_k \sum_{j=1}^{k-1} \no{U}_{H_x^s(\R)}^{j+1} \Big).
}

In view of \eqref{in1} and \eqref{in2}, and since $B_k \leq A_k$,  we infer
\eee{\label{in}
|I_{n}| 
\leq
c_k \Big(\no{U}_{H_x^s(\R)}^{k+2} + A_k \sum_{j=1}^{k-1} \no{U}_{H_x^s(\R)}^{j+2} \Big).
}
Combining \eqref{jn} and \eqref{in} with \eqref{injn}, we obtain
\eee{
\frac{d}{dt} \no{\p_x^nU}_{L_x^2(\R)}^2
\leq c_{n, k} \Big(\no{U}_{H_x^s(\R)}^{k+2} + A_k \sum_{j=1}^{k-1} \no{U}_{H_x^s(\R)}^{j+2} \Big)
}
so  summing over $0\le n\le s$ and integrating in $t$, we have
\eee{\label{int-est}
\no{U(t)}_{H_x^s(\R)}^2 
\leq
\no{U_0}_{H^s(\R)}^2
+
c_{s, k}   \int_0^t  \Big(\no{U(\tau)}_{H_x^s(\R)}^{k+2} + A_k \sum_{j=1}^{k-1} \no{U(\tau)}_{H_x^s(\R)}^{j+2} \Big) d\tau.
}

Next, along the lines of Section 3.1 in \cite{et2016book}, we let $T:=\inf\big\{t:\no{U(t)}_{H_x^s(\R)}\ge 2\no{U_0}_{H^s(\R)}\big\}$. Then, for $t\in [0, T]$ we have $\no{U(t)}_{H_x^s(\R)}\le 2\no{U_0}_{H^s(\R)}$ and so inequality \eqref{int-est} yields
\eee{\label{int-ineq2}
\no{U(t)}_{H_x^s(\R)}^2 
\leq
\no{U_0}_{H^s(\R)}^2
+
2^{k+2} c_{s, k}  \Big(\no{U_0}_{H^s(\R)}^{k+2} + A_k \sum_{j=1}^{k-1} \no{U_0}_{H^s(\R)}^{j+2} \Big) \, t.
}
Moreover, by continuity in time (since we are working with a smooth solution), it must be that $\no{U(T)}_{H_x^s(\R)}=2\no{U_0}_{H^s(\R)}$. Indeed, suppose $\no{U(T)}_{H_x^s(\R)}>2\no{U_0}_{H^s(\R)}$. Then, by continuity there exists an $\ve>0$ such that $\no{U(T-\ve)}_{H_x^s(\R)}>2\no{U_0}_{H^s(\R)}$, which contradicts the definition of $T$ as the infimum. Thus, evaluating \eqref{int-ineq2} at $t=T$ we obtain
\aaa{ 
4\no{U_0}_{H^s(\R)}^2 
&\leq
\no{U_0}_{H^s(\R)}^2
+
2^{k+2} c_{s, k} \Big(\no{U_0}_{H^s(\R)}^{k+2} + A_k \sum_{j=1}^{k-1} \no{U_0}_{H^s(\R)}^{j+2} \Big) \,  T
}
which can be rearranged to
\eee{\label{e:timebound}
T \geq
3 \, \Big[2^{k+2} c_{s, k}  \Big(\no{U_0}_{H^s(\R)}^k + A_k \sum_{j=1}^{k-1} \no{U_0}_{H^s(\R)}^j \Big) \Big]^{-1}.
}
Therefore, estimate \eqref{se} is valid until (at least) the lower bound on the right side of \eqref{e:timebound}. This concludes the argument  in the case of smooth solutions.

Next, we extend the result to the case of Sobolev solutions. To this end, by density, given $U_0 \in H^s(\R)$, we let $\{U_{0, n}\}_{n\in\N} \in \mathcal{S}(\R)$ be such that $\no{U_0-U_{0, n}}_{H^s(\R)}\to 0$ as $n\to\infty$. Then, thanks to the result of \cite{k1983}, there is a unique solution $U_n \in C_t^\infty([0, T_{0, n}]; \mathcal S(\R))$ corresponding to $U_{0, n}$, where $T_{0, n}>0$ depends only on the norm of $U_{0, n}$. Hence, by estimate \eqref{se}, which was established above for smooth solutions, 
\eee{
\no{U_n(t)}_{H_x^s(\R)} \leq 2 \no{U_{0, n}}_{H^s(\R)}, \quad 0 \leq t \leq T_n,
}
where  
$$
T_n := \frac{c_{s,k}}{\no{U_{0, n}}_{H^s(\R)}^k + A_k \sum_{j=1}^{k-1} \no{U_{0, n}}_{H^s(\R)}^j}.
$$
In turn, by the triangle inequality  and the Lipschitz continuity guaranteed by Theorem 3.2 of \cite{s1997}, we have
\begin{align*}
\no{U(t)}_{H_x^s(\R)} &= \no{U(t)-U_n(t)+U_n(t)}_{H_x^s(\R)}
\\
&\le \no{U_n(t)}_{H_x^s(\R)} + \no{U(t)-U_n(t)}_{H_x^s(\R)}
\\
&\le 2\no{U_{0, n}}_{H^s(\R)} + K\no{U_0-U_{0, n}}_{H^s(\R)},
\quad
t \in [0, T_n].
\end{align*}
Therefore, letting $n\to\infty$ yields the desired size estimate \eqref{se} valid for $0 \leq t \leq T$ where $T$ is given by \eqref{lspan}. This step completes the proof of Theorem \ref{se-t}.

\section{Proof of the distance estimates}
\label{s:proximity}

In this section, we establish Theorem \ref{prox-t} on the proximity between the solutions of different members of the gKdV family \eqref{gkdv} and the KdV equation \eqref{kdv}. Specifically, we estimate the  distance of solutions
\be
\label{e:diff}
\Delta(x,t) := u(x,t) - U(x,t)
\ee
in suitable function spaces, by considering the Cauchy problem satisfied by $\Delta$ and using the size estimate of Theorem~\ref{se-t} along the way. 
Subtracting  \eqref{gkdv} from  \eqref{kdv} yields the equation
\be 
\label{e:DeltaPDE}
\Delta_t + \Delta_{xxx} = N(x,t) := -uu_x + F(U)U_x. 
\ee
The Fourier transform of $f\in L_x^2(\R)$ is given by
\eee{
\widehat{f}(\xi) := \int_{\R}e^{-i\xi x}f(x)\,dx, \quad \xi\in\R,
}
with the usual inversion formula
\eee{
f(x):= \frac{1}{2\pi}\int_{\R}e^{i\xi x}\widehat{f}(\xi)\,d\xi, \quad x\in\R.  
}
Taking the Fourier transform of \eqref{e:DeltaPDE} and integrating the resulting ordinary differential equation in $t$, we obtain
\be
\label{ftd}
\widehat{\Delta}(\xi,t) = e^{i\xi^3t}\widehat{\Delta}(\xi,0) + \int_0^te^{i\xi^3(t-\tau)}\widehat{N}(\xi,\tau)\,d\tau.
\ee 
\noindent
\textbf{$L^2$ distance estimate.} We begin with the proof of the $L^2$ estimate \eqref{dl2}. By Plancherel's theorem, the triangle inequality and Minkowski's integral inequality, 
\eee{\label{dl2-0}
\no{\Delta(t)}_{L_x^2(\R)} 
\leq
\no{\Delta(0)}_{L^2(\R)} 
+ 
\int_0^t \no{N(\tau)}_{L_x^2(\R)} d\tau.
}
Next, we estimate $N$ in $L_x^2(\R)$. By the triangle inequality and  the Sobolev embedding theorem,
\aaa{
\no{N}_{L_x^2(\R)} 
&\leq
\no{uu_x}_{L_x^2(\R)} + \no{F(U)U_x}_{L_x^2(\R)}
\nn\\
&\leq
\no{u}_{L_x^{\infty}(\R)} \no{u_x}_{L_x^2(\R)} 
+ 
\sum_{j=1}^k |a_j| \no{U^j}_{L_x^\infty(\R)} \no{U_x}_{L_x^2(\R)}
\nn\\
&\leq
\no{u}_{H_x^1(\R)}^2
+ 
\sum_{j=1}^k |a_j| \no{U^j}_{H_x^1(\R)} \no{U}_{H_x^1(\R)}.
}
Then, by the algebra property in $H^1(\R)$, 
\eee{
\no{N}_{L_x^2(\R)} 
\leq
\no{u}_{H_x^1(\R)}^2
+ 
c \sum_{j=1}^k |a_j| \no{U}_{H_x^1(\R)}^{j+1}
\leq
\no{u}_{H_x^1(\R)}^2
+ 
c \, \Big(|a_k| \no{U}_{H_x^1(\R)}^{k+1} + B_k \sum_{j=1}^{k-1}  \no{U}_{H_x^1(\R)}^{j+1} \Big)
}
with $B_k$ given by \eqref{bkdef}. In turn, recalling the definition of $c_k$ in \eqref{jn}, 
\eee{\label{e:NL2}
\int_0^t\no{N(\tau)}_{L_x^2(\R)} d\tau
\leq 
\sup_{\tau \in [0, t]}
\Big[
\no{u(\tau)}_{H_x^1(\R)}^2
+ 
c_k  \Big(\no{U(\tau)}_{H_x^1(\R)}^{k+1} + B_k \sum_{j=1}^{k-1}  \no{U(\tau)}_{H_x^1(\R)}^{j+1} \Big)  \Big] \, t
}
thus, returning to \eqref{dl2-0}, we infer  
\eee{
\label{dhs1}
\no{\Delta(t)}_{L_x^2(\R)} 
\leq 
\no{\Delta(0)}_{L^2(\R)} 
+ 
\sup_{\tau \in [0, t]}
\Big[
\no{u(\tau)}_{H_x^1(\R)}^2
+ 
c_k  \Big(\no{U(\tau)}_{H_x^1(\R)}^{k+1} + B_k \sum_{j=1}^{k-1}  \no{U(\tau)}_{H_x^1(\R)}^{j+1} \Big)  \Big] \, t.
}
Replacing the $H_x^1(\R)$ norms on the right side by $H_x^2(\R)$ norms and then taking the supremum over $t\in [0, T]$ for any fixed $T \in (0, T_1]$ where $T_1$ is given by \eqref{tcs}, we invoke the size estimates \eqref{se} and \eqref{sec} for the gKdV and KdV solutions (this is possible because $T\leq T_1$ and $T_1$ is the minimum of \eqref{lspan} and \eqref{lspank} for $s=2$ and, in the latter case, $k=1$) to obtain
\eee{
\label{e:keyest1T}
\sup_{t\in [0, T]} \no{\Delta(t)}_{L_x^2(\R)} 
\leq 
\no{\Delta(0)}_{L^2(\R)} 
+ 
2^{k+1} \Big[\no{u_0}^2_{H^2(\R)} + c_k  \Big(\no{U_0}_{H^2(\R)}^{k+1} + B_k \sum_{j=1}^{k-1}  \no{U_0}_{H^2(\R)}^{j+1} \Big) \Big] \, T.
}
Hence, in view of \eqref{icorder} and \eqref{icdorder}, we deduce 
\eee{
\label{IN1}
\sup_{t\in [0, T]} \no{\Delta(t)}_{L_x^2(\R)} 
\leq 
C \max\left\{\epsilon^2, \epsilon^{k+1}\right\} 
+ 
2^{k+1} \Big[c_0^2 \epsilon^2 +  c_k  \Big(C_0^{k+1} \epsilon^{k+1} + B_k \sum_{j=1}^{k-1}  C_0^{j+1} \epsilon^{j+1} \Big)\Big] \, T
}
which readily yields the claimed distance estimate \eqref{dl2}.
\\[2mm]
\textbf{$H^s$ distance estimate for any $s\in\N$.}
We proceed to the proof of the $H^s$ estimate \eqref{dhs}. 
By the Fourier transform definition of the Sobolev norm, the triangle inequality and Minkowski's integral inequality, 
\eee{\label{dhs-0}
\no{\Delta(t)}_{H_x^s(\R)} 
\leq
\sqrt 2 
\no{\Delta(0)}_{H^s(\R)} 
+ 
\sqrt 2 
\int_0^t \no{N(\tau)}_{H_x^s(\R)} d\tau.
}
Since $s>\frac 12$, we apply the algebra property to obtain
\aaa{
\no{N}_{H_x^s(\R)}
&\leq
\no{uu_x}_{H_x^s(\R)} + \no{F(U) U_x}_{H_x^s(\R)}
\nn\\
&\leq
c_s \no{u}_{H_x^s(\R)} \no{u}_{H_x^{s+1}(\R)}
+
c_s \no{F(U)}_{H_x^s(\R)} \no{U}_{H_x^{s+1}(\R)}.
}
Furthermore, once again by the triangle inequality and the algebra property,
\eee{
\no{F(U)}_{H_x^s(\R)}
\leq
\sum_{j=1}^k |a_j| \no{U^j}_{H_x^s(\R)}
\leq
c_s \sum_{j=1}^k |a_j| \no{U}_{H_x^s(\R)}^j
\leq
c_s c_k \Big(\no{U}_{H_x^s(\R)}^k + B_k  \sum_{j=1}^{k-1} \no{U}_{H_x^s(\R)}^j \Big)
}
with $c_k$ and $B_k$ defined by \eqref{jn} and \eqref{bkdef}. 
Therefore, 
\aaa{
\no{N}_{H_x^s(\R)}
&\leq
c_{s, k}
\Big[
\no{u}_{H_x^s(\R)} \no{u}_{H_x^{s+1}(\R)}
+
\no{U}_{H_x^{s+1}(\R)} \Big(\no{U}_{H_x^s(\R)}^k + B_k  \sum_{j=1}^{k-1} \no{U}_{H_x^s(\R)}^j \Big) 
\Big]
\nn\\
&\leq
c_{s, k} 
\Big(
\no{u}_{H_x^{s+1}(\R)}^2
+
\no{U}_{H_x^{s+1}(\R)}^{k+1} + B_k  \sum_{j=1}^{k-1} \no{U}_{H_x^{s+1}(\R)}^{j+1} 
\Big).
}
Note that we have boosted some of the Sobolev exponents on the right side  in order to have at least $H_x^2(\R)$ (since we are working with $s \in \N$) and thus be able to employ Theorem \ref{se-t}. In turn, 
\eee{
\label{lgc}
\no{\Delta(t)}_{H_x^s(\R)}
\leq 
\sqrt 2 \no{\Delta(0)}_{H^s(\R)}
+
\sqrt 2  
\sup_{\tau \in [0, t]}
c_{s, k} 
\Big(
\no{u(\tau)}_{H_x^{s+1}(\R)}^2
+
\no{U(\tau)}_{H_x^{s+1}(\R)}^{k+1} + B_k  \sum_{j=1}^{k-1} \no{U(\tau)}_{H_x^{s+1}(\R)}^{j+1} 
\Big) \, t
}
and taking the supremum over $t\in [0, T]$ for any fixed $T \in (0, T_s]$ with $T_s$ given by \eqref{tcs}, we invoke the size estimates~\eqref{se} and \eqref{sec} for the gKdV and KdV solutions (this is possible because $T\leq T_s$ and $T_s$ is the minimum of \eqref{lspan} and \eqref{lspank} with $k=1$ in the latter case) to infer
\aaa{
\sup_{t\in [0, T]} \no{\Delta(t)}_{H_x^s(\R)}
&\leq
\sqrt 2 \no{\Delta(0)}_{H^s(\R)}
+
 2^{k+\frac 32} 
 \Big(
\no{u_0}_{H^{s+1}(\R)}^2
+
\no{U_0}_{H^{s+1}(\R)}^{k+1} + B_k  \sum_{j=1}^{k-1} \no{U_0}_{H^{s+1}(\R)}^{j+1} 
\Big) \, T.
}
Therefore, employing \eqref{icorders} and \eqref{icdorders}, we conclude that
\eee{
\label{IN2}
\sup_{t\in [0, T]} \no{\Delta(t)}_{H_x^s(\R)} 
\leq 
\sqrt 2 \, C_s \max\left\{\epsilon^2, \epsilon^{k+1}\right\} 
+ 
2^{k+\frac 32}  \Big( c_{0, s}^2 \epsilon^2 + C_{0, s}^{k+1} \epsilon^{k+1} + 
B_k  \sum_{j=1}^{k-1} C_{0, s}^{j+1} \epsilon^{j+1} \Big) \, T
}
which implies the claimed distance estimate \eqref{dhs},  completing the proof of Theorem \ref{prox-t}.

\section{Numerical simulations}
\label{s:numerical}

The numerical experiments are devoted to two important members of the gKdV family \eqref{gkdv}: 
(i) the standard \textit{monomial} case
\be
\label{gkdvp}
F(U)= U^k, \quad k\in\mathbb{N}, 
\ee 
with the special cases $k=1$ and $k=2$ corresponding to the KdV and mKdV equations, which are the only integrable members of the family \eqref{gkdv}, and (ii) the \textit{weakly perturbed} cases
\be
\label{wkdv}
F(U)=U+\delta U^k \text{ \ and \ } F(U)=U^2+\delta U^m, \quad k,  m\in\mathbb{N}, \ 0<\delta<1, 
\ee
which correspond to  weakly non-integrable KdV and mKdV equations, respectively, in the spirit of \cite{dz2002} (recall that $k=2$ and $\delta>0$ corresponds to the integrable Gardner equation).

In the case of the monomial power nonlinearity \eqref{gkdvp}, the gKdV family \eqref{gkdv} admits the one-soliton solution
\be
\label{e:soliton}
U(x,t) = \left(\frac{c(k+1)(k+2)}{2}\right)^{\frac 1k}\sech^{2k}\left(\frac{k\sqrt{c}}{2}(x-ct-x_0)\right),
\ee
where $c>0$ and $x_0$ respectively denote the wave speed and  the initial position (or phase). In the case of the integrable KdV ($k=1$) and mKdV ($k=2$) equations, the one-soliton \eqref{e:soliton} yields two useful candidates for the initial data of our simulations, namely
\eee{\label{e:kdvsoliton}
u_{\text s}(x) := 3c\sech^2\left(\frac{\sqrt{c}}{2}(x-x_0)\right), 
\quad
U_{\text s}(x) := \sqrt{6c}\sech\left(\sqrt{c}\left(x-x_0\right)\right),
} 
where we shall take $x_0=0$ without loss of generality in the case of one-soliton solutions. The Sobolev $H^s$ norm with $s=0,1,2$ can be computed in closed form for both of the above initial profiles. In particular, 
\be
\label{e:kdvH2}
\no{u_{\text s}}_{H^2}^2 = 24\sqrt{c} \, \Big(c + \frac{1}{5}c^2 + \frac{1}{7}c^{3}\Big),
\quad
\no{U_{\text s}}_{H^2}^2 = 4\sqrt{c}\, \Big(3 + c + \frac{7}{5}c^{2}\Big).
\ee
This fact will be useful for the numerical demonstration of our proximity results.

The numerical method used in the implementation of our experiments is a spectral approach with Fourier-type discretizations on the interval $x\in[-L,L]$ for $L$ sufficiently large. Time-stepping is handled using a fourth order Runge-Kutta method (RK4).
\begin{figure}[t!]
	\centering
	\begin{tabular}{cc}
		\textbf{(a)}\hspace{5.0cm}\textbf{(b)}\\	
		\includegraphics[width=6cm]{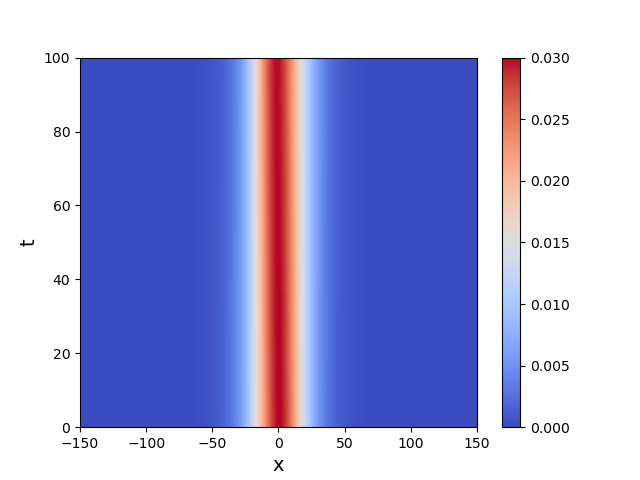}
		\includegraphics[width=6cm]{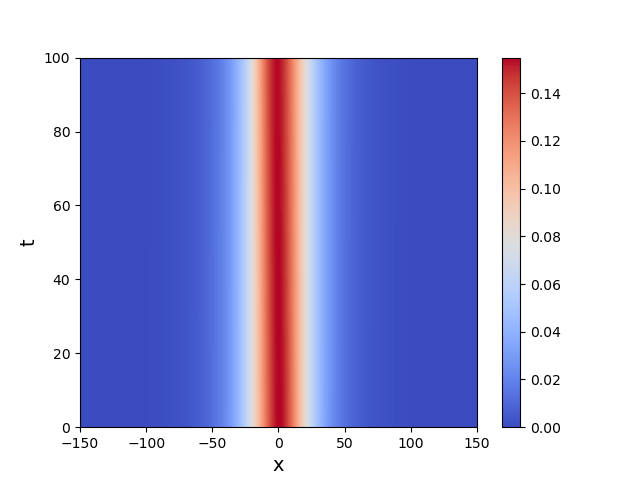}
		\\
		\textbf{(c)}\hspace{5.0cm}\textbf{(d)}\\	
		\includegraphics[width=6cm]{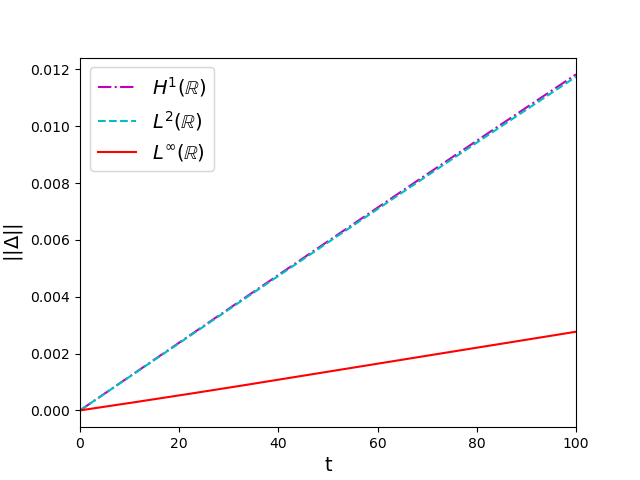}
		\includegraphics[width=6cm]{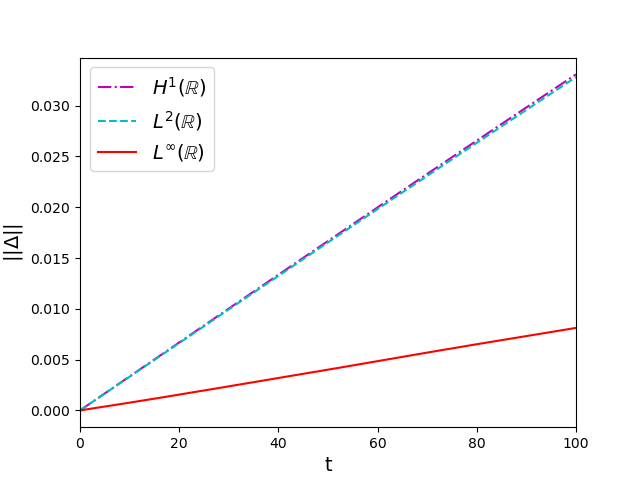}
	\end{tabular}
	\caption{The solution $U(x, t)$ of the gKdV equation \eqref{gkdv} with  monomial nonlinearity according to \eqref{gkdvp}. (a) The supercritical case $k=5$ supplemented with KdV one-soliton initial data $u_{\text s}(x)$ with $c=0.01$ and $x_0=0$. (b) The critical case $k=4$ supplemented with mKdV one-soliton initial data $U_{\text s}(x)$ with $c=0.004$ and $x_0=0$.
		(c) Time evolution of the norms $\no{\Delta(t)}_{\mathcal{X}}$ of the difference between the above gKdV solutions and the corresponding KdV and mKdV one-solitons for $\mathcal{X}= H_x^1$, $L_x^2$ and $L_x^{\infty}$, for $k=5$. (d) The corresponding time evolution of the norms for $k=4$.
	}
	\label{f:1}
\end{figure}
\subsection{The gKdV with power nonlinearity}\label{41-ss}

 We begin the illustration of our proximity results stated in Theorem~\ref{prox-t} and Corollary~\ref{prox-c} by first considering the case of the gKdV equation with power nonlinearity \eqref{gkdvp}. We examine the dynamics of the two initial conditions $u_{\text s}(x)$ and $U_{\text s}(x)$ in \eqref{e:kdvsoliton}, which correspond to the one-soliton solutions of the KdV and mKdV equations, respectively.
\\[2mm]
\noindent
\textbf{I. Small initial data.}
We begin with the initial data \eqref{e:kdvsoliton} in the case of small size $\no{u_{\text s}}_{H^2} = \no{U_{\text s}}_{H^2} = \epsilon<1$. These initial data are evolved numerically according to the non-integrable gKdV equation with power nonlinearity~\eqref{kdvk}, and the resulting dynamics are compared with those of the same initial conditions evolved by the integrable KdV equation. For the gKdV case, we focus on the critical exponent $k=4$ and the supercritical exponent $k=5$.

As a first example, we consider the initial datum $u_{\text s}(x)$ in \eqref{e:kdvsoliton} with $c=0.01$ so that, according to \eqref{e:kdvH2}, $\no{u_{\text s}}_{H^2} \approx 0.1551$. We then evolve this datum under the supercritical gKdV equation \eqref{kdvk} with $k=5$. The dynamics for this case are depicted in Figure~\ref{f:1}, where panel (a) illustrates the contour plot of the solution $U(x,t)$ of the gKdV equation and  panel (c) illustrates the time evolution of the norms $\no{\Delta(t)}_{\mathcal{X}}$.  Evolving this datum over the time interval $t\in [0, \no{u_{\text s}}_{H^2}^{-1}] \simeq [0, 6.5]$ and  computing numerically the difference norm $\no{\Delta(t)}_{\mathcal{X}}$ for $\mathcal{X}=H_x^s(\mathbb{R})$, $s=0,1$, and  $L_x^{\infty}(\mathbb{R})$, we observe the following: By Corollary~\ref{prox-c}, we have the theoretical bound $\no{\Delta(t)}_{\mathcal{X}} \lesssim \epsilon^2 \approx 0.024$. The numerical simulations fulfill this proximity estimate as shown in panel (c) of  Figure~\ref{f:1} and, in fact, for $t \in [0, 80]$ yield $\no{\Delta(t)}_{L_x^{\infty}}\le 0.0024$, which is an order of magnitude smaller than the theoretical bound. Moreover, the simulations exhibit a linear in time growth of the proximity error, in  agreement with the theoretical predictions.

\begin{figure}[t!]
\centering
\begin{tabular}{cc}
	\textbf{(a)}\hspace{5.0cm}\textbf{(b)}\hspace{5.0cm}\textbf{(c)}\\	
\includegraphics[width=5cm]{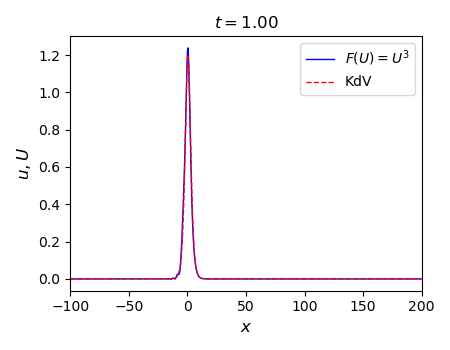}
\includegraphics[width=5cm]{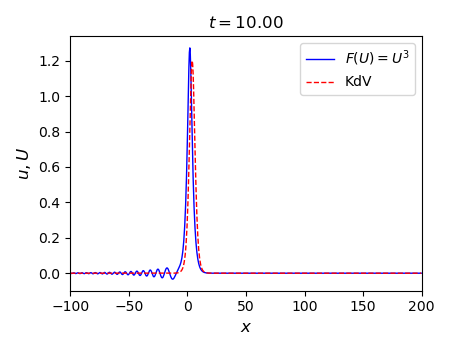}
\includegraphics[width=6cm]{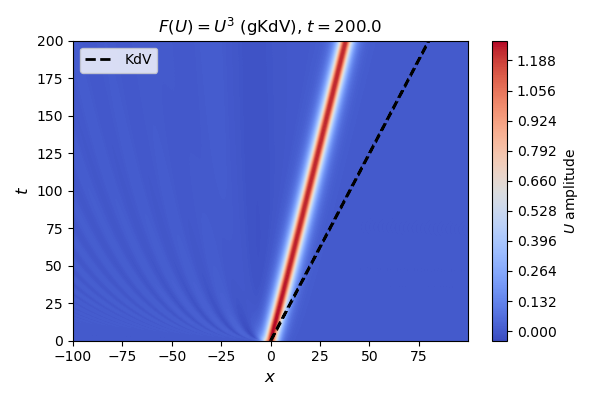}\\
\textbf{(d)}\hspace{5.0cm}\textbf{(e)}\hspace{5.0cm}\textbf{(f)}\\	
\includegraphics[width=5cm]{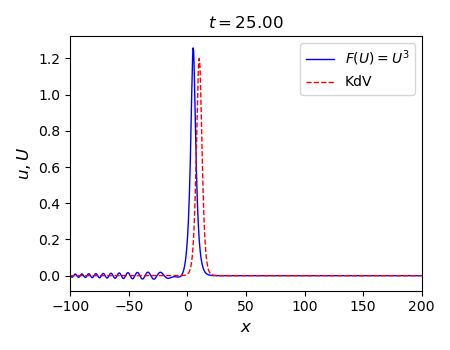}
\includegraphics[width=5cm]{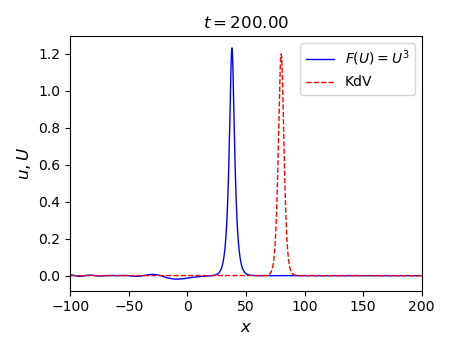}
\includegraphics[width=4.7cm]{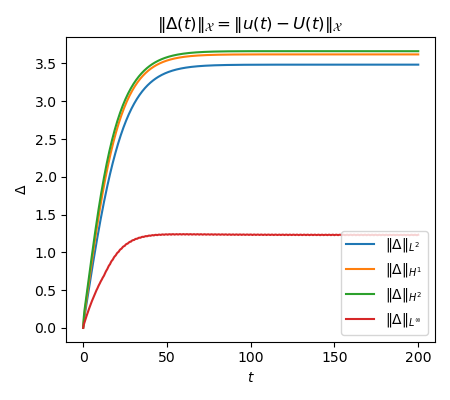}
\hspace{1cm}
\end{tabular}
\caption{(a), (b), (d) and (e) Snapshots comparing the profile of the solution of gKdV with power nonlinearity \eqref{kdvk} for $k=3$ (solid blue) against the profile of the solution of the  integrable KdV~\eqref{kdv} (red dashed). (c) Contour plot of the evolution of the initial  datum $u_{\text s}(x)$ in \eqref{e:kdvsoliton} with $c=0.4$ and $x_0=0$ under gKdV. The dashed black line is the caustic of the analytical KdV soliton. (f) The evolution of the distance norms in $\mathcal{X}= H_x^s$, $s=0,1,2$, and $L_x^{\infty}$.
}
\label{fig4.2}
\end{figure}

As a second example, we consider the second initial profile $U_{\text s}(x)$  in \eqref{e:kdvsoliton} with $c=0.004$ so that, by \eqref{e:kdvH2},  $\no{U_{\text s}}_{H^2} \approx 0.8718$. We evolve the mKdV one-soliton data under the critical gKdV equation \eqref{kdvk} with $k=4$. The corresponding dynamics are depicted in panels (b) and (d) of Figure~\ref{f:1}; the format of the presentation is the same as in panels (a) and (c), respectively. For $t\in [0,\no{U_{\text s}}_{H^2}^{-2}] \approx [0, 1.32]$, we compute numerically $\no{\Delta(t)}_{\mathcal{X}}$. By Corollary~\ref{prox-c}, we have $\no{\Delta(t)}_{\mathcal{X}} \lesssim \epsilon^2 \approx 0.76$. Once again, the numerical results fulfill the theoretical estimates. Specifically, we find $\no{\Delta(t)}_{L_x^{\infty}}\le 0.01$ for $t\in [0, 100]$, which is again about an order of magnitude smaller than the theoretical bound. As in the first example, the growth of the distances in panel (d) of Figure~\ref{f:1} is in agreement with the theoretically predicted linear temporal growth of the proximity estimates.
\\[2mm]
\noindent
\textbf{II. Large initial data: the effect of rotation of the linear caustic of the KdV soliton.}  For solitons of higher amplitude such that $\no{u_{\text s}}_{H^2} = \no{U_{\text s}}_{H^2} = \epsilon>1$, according to Corollary \ref{prox-c} we expect a rapid deviation between the dynamics of the integrable and the non-integrable equations. This deviation is depicted in Figure~\ref{fig4.2}, which compares the dynamics between the gKdV equation \eqref{kdvk} with $k=3$ and the integrable KdV equation for the one-soliton initial datum $u_{\text s}(x)$ in $\eqref{e:kdvsoliton}$ with $c=0.4$, for which $\no{u_{\text s}}_{H^2}=\epsilon\approx 2.58$. The system is integrated for $x\in [-L, L]$ with $L=1000$ and $t\in [0, 4000]$.

Panels (a), (b), (d) and (e) in Figure~\ref{fig4.2} illustrate snapshots at specific times comparing the profiles of the solutions. The gKdV solution is shown with the solid blue curve and the KdV soliton with the dashed red curve. Panel (c) shows the contour plot of the gKdV solution, together with the caustic of the corresponding analytical KdV soliton, depicted by the dashed black line. Panel (f) shows the evolution of the distance norm $\no{\Delta(t)}_{\mathcal{X}}$ with $\mathcal{X}=H_x^s(\mathbb{R})$, $s=0,1,2$, and $L_x^{\infty}(\mathbb{R})$.

According to Corollary \ref{prox-c}, we expect proximity of the solutions only for very small timescales predicted by $\widetilde T_s=\mathcal{O}(\no{u_{\text s}}_{H^2}^{-3})$. This is already seen in the first snapshot at $t=1$ of panel (a), where we observe the beginning of the deviation of the numerical gKdV soliton  from the analytical KdV one-soliton, manifested by the small increase in its amplitude and the onset of its radiation. In the snapshot at $t=10$ of panel (b), the radiation of the gKdV soliton is clearly observable and the  KdV soliton is faster.
These snapshots, as well as the one at $t=25$ of panel (d), depict profiles in the regime of the linear growth of the distances between the solutions. In particular, the snapshot at $t=25$ refers to a time close to the end of the linear growth regime of $\no{\Delta(t)}_{\mathcal{X}}$ shown in panel (f), manifested by the complete separation of the soliton profiles, as observed in the snapshot at $t=200$ of panel (e); the beginning of this separation at $t\sim 25$ is also visible in the contour plots of the solutions shown in panel (c).

We observe that, apart from the radiative effect (which, although of small amplitude, contributes to the linear growth of the distance in the Sobolev norms), the core of the gKdV solitary wave maintains a close amplitude and profile to the analytical KdV soliton, though the gKdV wave is noticeably slower. This difference in speed leads to the different angles of the caustics of the gKdV and KdV solutions, visually confirming the velocity difference between the faster KdV soliton and the slower gKdV solitary wave. The norms $\no{\Delta(t)}_{\mathcal{X}}$ achieve a constant value after the separation, as they essentially become the sum of the individual norms of the gKdV and KdV solitons (in the case of the $L^\infty$ norm, only the largest norm of the gKdV soliton is measured).

This observation motivated us to implement a ``left rotation'' of the linear caustic of the integrable KdV soliton in order to reduce the speed of the analytical soliton so that its caustic will coincide with that of the gKdV solitary wave. Indeed, this can be done by considering the \textit{rescaled} integrable KdV equation 
\be
\label{resckdv}
u_t + \nu u_{xxx} + \nu u u_{x} = 0
\ee
with the corresponding rescaled analytical one-soliton solution
\be
\label{e:rsoliton}
u_{s,\nu}(x,t) = 3c\,\sech^{2}\left(\frac{\sqrt{c}}{2}\,(x - \nu c t - x_0)\right),
\ee
obtained by applying the transformation $t \mapsto t/\nu$ to the standard KdV equation~\eqref{kdv}.

\begin{figure}[t!]
\centering
\begin{tabular}{cc}
\textbf{(a)}\hspace{5.0cm}\textbf{(b)}\hspace{5.0cm}\textbf{(c)}\\	
\includegraphics[width=5cm]{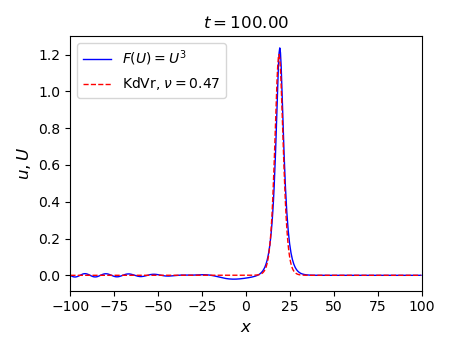}
\includegraphics[width=5cm]{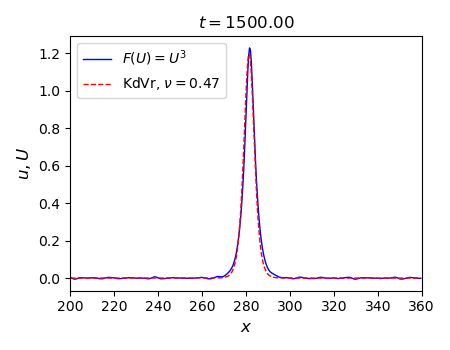}
\includegraphics[width=6cm]{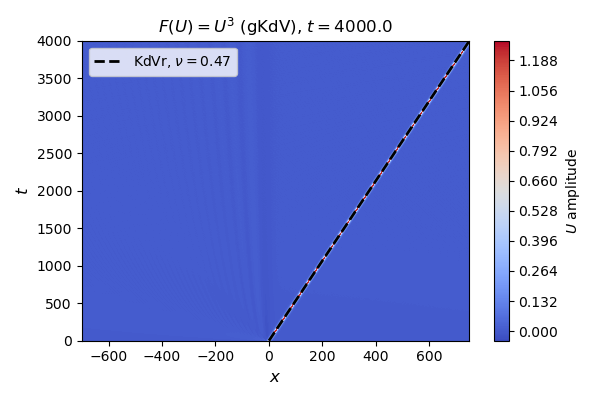}\\
\textbf{(d)}\hspace{5.0cm}\textbf{(e)}\hspace{5.0cm}\textbf{(f)}\\	
\hspace{-1.1cm}\includegraphics[width=5cm]{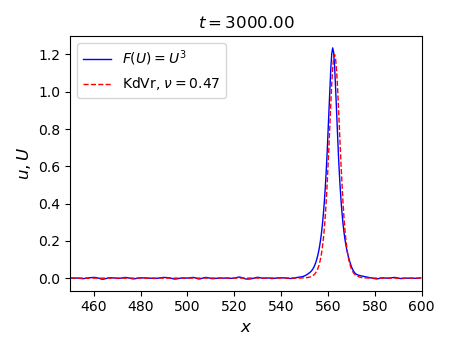}
\includegraphics[width=5cm]{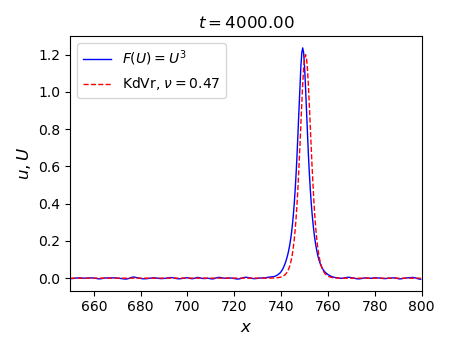}
\includegraphics[width=4.7cm]{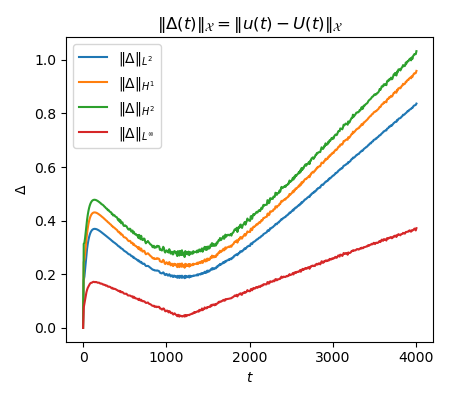}
\end{tabular}
\caption{ (a), (b), (d) and (e) Snapshots comparing the profiles of the solution of the gKdV equation with power nonlinearity~\eqref{kdvk} for $k=3$ (solid blue) against the profiles of the soliton~\eqref{e:rsoliton} of the  rescaled integrable KdV \eqref{resckdv} with $\nu=0.47$ (red dashed). (c) Contour plot of the evolution of the initial datum \eqref{e:kdvsoliton} with $c=0.4$ and $x_0=0$ under gKdV. The dashed black line is the caustic of the rescaled soliton \eqref{e:rsoliton} with $\nu=0.47$. (f) The evolution of the distance measured in the spaces $\mathcal{X}= H_x^s$, $s=0,1,2$, and $L_x^{\infty}$. 
}
\label{fig4.3}
\end{figure}

\begin{figure}[t!]
\centering
\begin{tabular}{cc}
\textbf{(a)}\hspace{5.0cm}\textbf{(b)}\hspace{5.0cm}\textbf{(c)}\\	
\includegraphics[width=5cm]{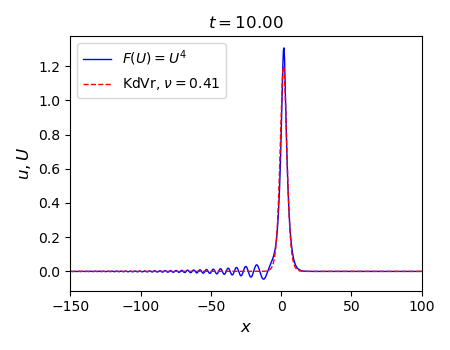}
\includegraphics[width=5cm]{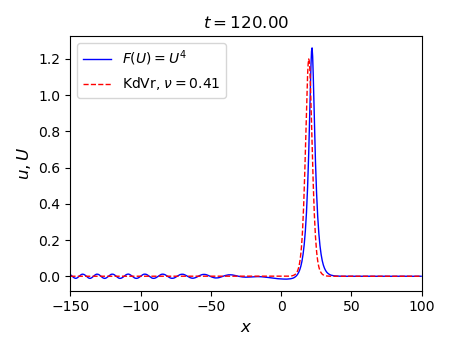}
\includegraphics[width=5.5cm]{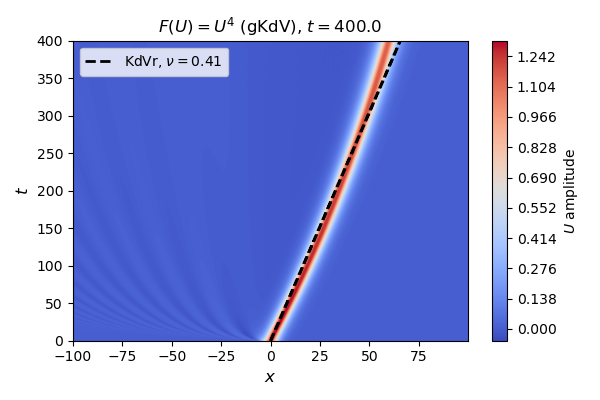}\\
\textbf{(d)}\hspace{5.0cm}\textbf{(e)}\hspace{5.0cm}\textbf{(f)}\\	
\includegraphics[width=5cm]{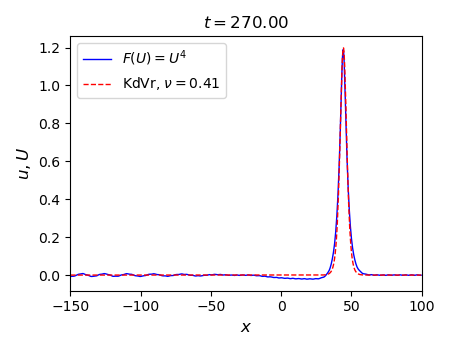}
\includegraphics[width=5cm]{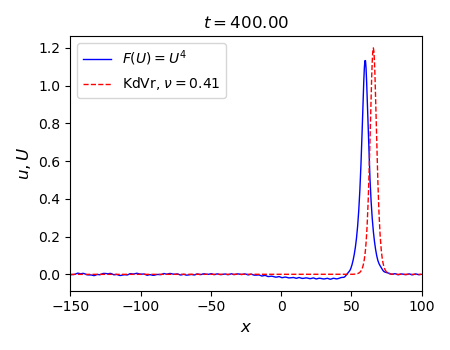}
\includegraphics[width=4.7cm]{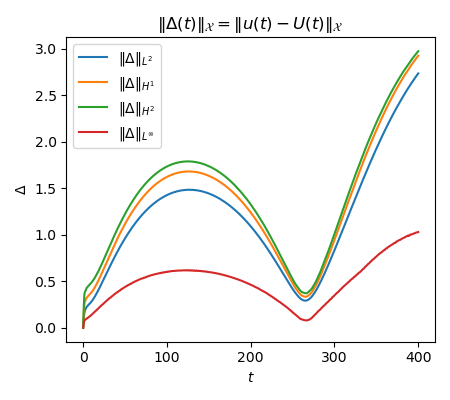}\hspace{0.7cm}
\end{tabular}
\caption{(a), (b), (d) and (e) Snapshots comparing the profiles of the gKdV solution in the case of the power nonlinearity \eqref{kdvk} with $k=4$ (solid blue) against the profiles of the soliton \eqref{e:rsoliton} of the  rescaled integrable KdV \eqref{resckdv} with $\nu=0.41$ (red dashed). (c) Contour plot of the evolution of the initial condition \eqref{e:kdvsoliton} with $c=0.4$, $x_0=0$ under gKdV. The dashed black line is the caustic of the rescaled soliton \eqref{e:rsoliton} with $\nu=0.41$. (f) The evolution of norms of the distance between the solutions $\mathcal{X}= H_x^s$, $s=0,1,2$, and $L_x^{\infty}$.  
}
\label{fig4.4}
\end{figure}

We detect the following dynamical feature: for a given nonlinearity exponent $k \geq 2$, there exists a scaling parameter $0<\nu(k) < 1$ such that the soliton dynamics in the gKdV equation remain remarkably proximal, for significantly large time intervals, to the analytical rescaled soliton~\eqref{e:rsoliton}.  
For gKdV with $k=3$, we detect a scaling parameter $\nu(3)\sim0.47$. Figure~\ref{fig4.3} depicts the comparison of the dynamics between the gKdV and the rescaled KdV equations for the initial datum obtained from \eqref{e:rsoliton} with $c=0.4$ and $x_0=0$; the format of the presentation is the same as in Figure~\ref{fig4.2}. We observe in panel (f) of Figure~\ref{fig4.3} that the timescales required for the distance $\no{\Delta(t)}_{H_x^s(\mathbb{R})}$ to reach $\mathcal{O}(1)$ are of $\mathcal{O}(10^3)$ ($t\sim 4000$ in this example). In the critical case $k=4$, for which the relevant dynamics are shown in Figure~\ref{fig4.4}, we detect that $\nu(4)\sim 0.41$  and the timescales for the distance to reach $\mathcal{O}(1)$  are $t \sim \mathcal{O}(10^2)$.  
Within these large timescales (in view of the degree of the nonlinearity as discussed above), the observed deviation is primarily due to weak radiative effects exhibited by the gKdV dynamics, while the soliton’s initial amplitude and shape persist with only minor differences from the analytical rescaled soliton profile.

Inspecting the evolution of the distance norms  in Figures \ref{fig4.2} and \ref{fig4.4} for both the $k=3$ and $k=4$ examples, we observe a distinct dynamical pattern within the time interval of integration. The distance first reaches a local  maximum at the end of the initial linear growth regime. Subsequently, the distance decreases, attaining a minimum, after which the distance resumes an almost linear growth. This feature can be explained as follows. The power nonlinearity affects the KdV soliton initial datum by introducing ``curvature'' to its caustic, along which the solitary wave moves at a significantly lower speed than the analytical KdV soliton (as seen in panel (c) of Figure \ref{fig4.2} for $k=3$, prior to the caustic rotation). This feature is more clearly visualized when comparing the contour plot and evolution of norms as seen in panels (c) and (f) of Figure \ref{fig4.4} for $k=4$, where the curvature of the gKdV caustic is clearly visible. For the case $k=3$, the curvature of the gKdV caustic is very small relative to the time of evolution, causing the caustics to appear as almost coinciding.
After the rotation, the caustic of the rescaled KdV soliton (dashed black lines in the relevant contour plots) acts as a secant line traced from $t=0$, which approximates the slightly curved caustics of the gKdV solitary waves over large timescales. We relate the evolution of the norms to the relative speeds of the caustics:
\begin{enumerate}[label=$\bullet$, leftmargin=3.5mm, itemsep=2mm, topsep=1mm]
    \item \underline{Linear growth regime (1)}: When the dashed black  secant line is slightly \textit{above} the gKdV caustic, the gKdV soliton evolves with slightly higher speed and amplitude. This disparity defines the initial linear growth of the norms.
    \item \underline{Maximum and decrease}: The maximum of the norms refers to the maximum distance between the caustics. Following this maximum, the distance decreases as the secant line tends toward the crossing point with the curved gKdV caustic. At this crossing, both solitons have momentarily attained the same speed and amplitude, resulting in the minimum distance norm.
    \item \underline{Linear growth regime (2)}: Beyond this minimum, the linear secant caustic departs from the curved caustic of the gKdV soliton to the right, and the regime of linear growth of the distances begins again. In this phase, the rescaled KdV soliton becomes faster and gains amplitude relative to the gKdV solitary wave.
\end{enumerate}
\subsection{The weakly perturbed KdV and mKdV equations}

\begin{figure}[t!]
\centering
\begin{tabular}{cc}
\textbf{(a)}\hspace{5.0cm}\textbf{(b)}\hspace{5.0cm}\textbf{(c)}\\	
\includegraphics[width=5cm]{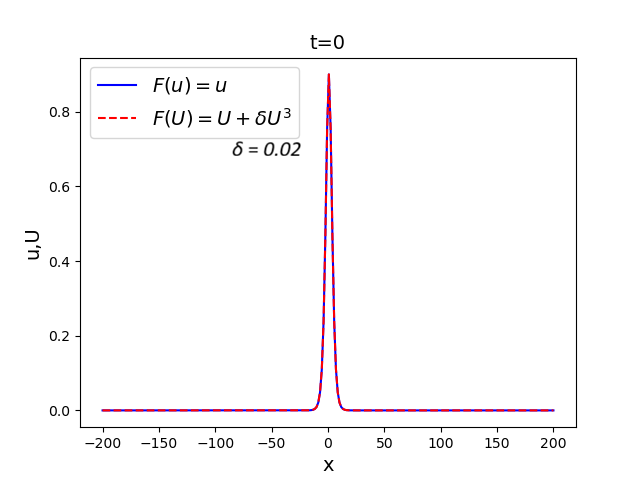}
\includegraphics[width=5cm]{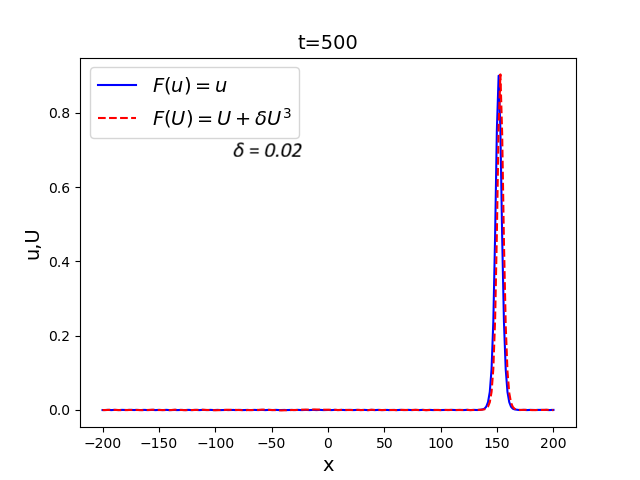}
\includegraphics[width=5cm]{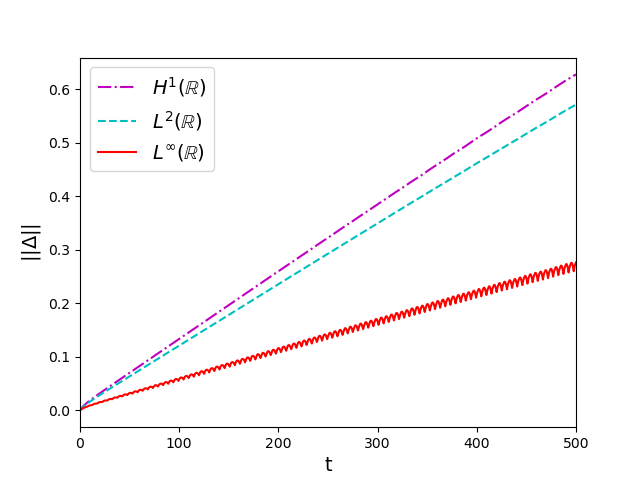}
\\
\textbf{(d)}\hspace{5.0cm}\textbf{(e)}\hspace{5.0cm}\textbf{(f)}\\
\includegraphics[width=5cm]{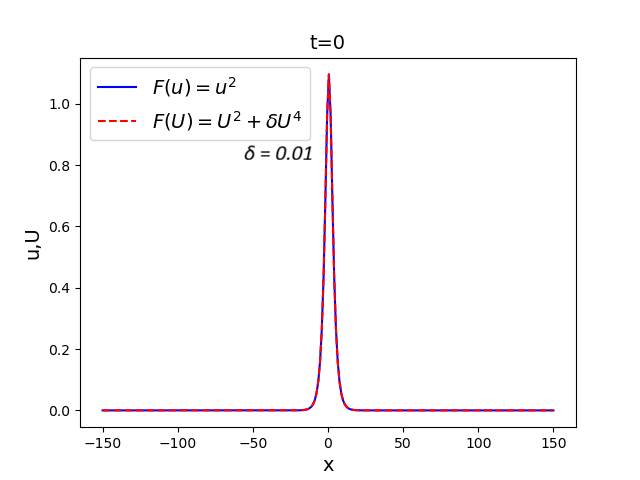}
\includegraphics[width=5cm]{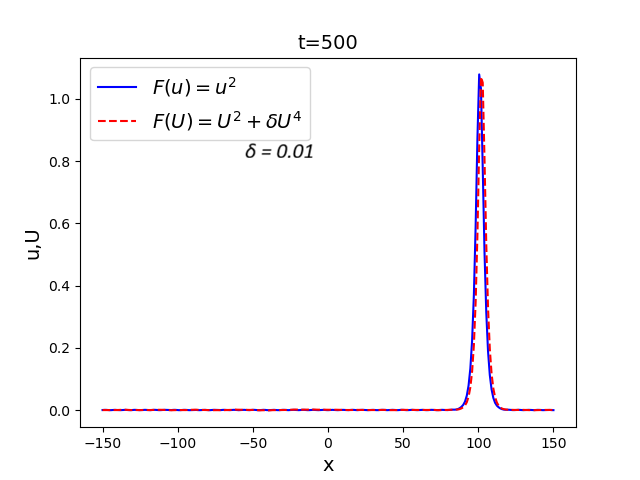}
\includegraphics[width=5cm]{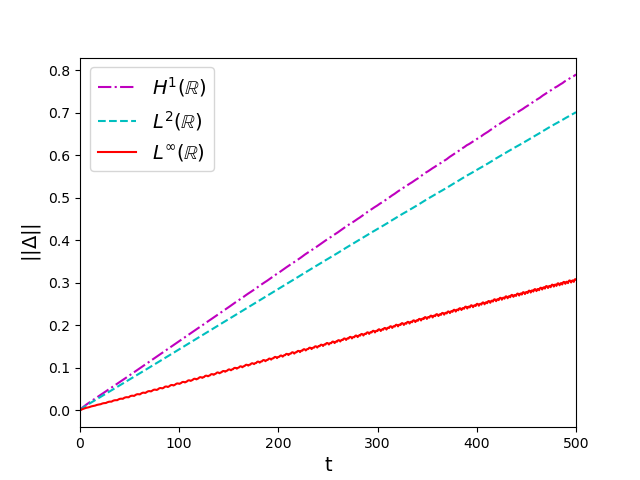}
\end{tabular}
\caption{(a) and (b) Snapshots comparing the profile of the solution of the gKdV equation \eqref{gkdv} with $F(U)=U+\delta U^3$ and $\delta=0.02$, supplemented with the KdV one-soliton initial datum $u_{\text s}(x)$ in \eqref{e:kdvsoliton} with $c=0.3$ and $x_0=0$, against the profile of the corresponding analytical KdV  one-soliton. Panel (a) shows the initial profiles and panel (b) the solutions at $t=500$. (c) The corresponding evolution of the norms $\no{\Delta(t)}_{\mathcal{X}}$ with $\mathcal{X}=H_x^1$, $L_x^2$ and $L_x^{\infty}$. (d) and (e)  Snapshots comparing the profile of the solution of the gKdV  \eqref{gkdv} with $F(U)=U^2+\delta U^4$ and$\delta=0.01$, supplemented with the mKdV one-soliton initial datum $U_{\text s}(x)$ in \eqref{e:kdvsoliton} with $c=0.2$ and $x_0=0$, against the profile of the corresponding analytical mKdV one-soliton. Panel (d) shows the initial profiles and panel (e) the solutions at $t=500$.  (f) The corresponding time evolution of the norms $\no{\Delta(t)}_{\mathcal{X}}$ with $\mathcal{X}=H_x^1$, $L_x^2$ and $L_x^{\infty}$.
}
\label{f:2}
\end{figure}

We now consider the weakly perturbed cases of the KdV and mKdV equations, where the gKdV equation is defined according to \eqref{wkdv}. For the numerical experiments in this setting, we consider both one-soliton and two-soliton initial data.
\\[2mm]
\noindent
\textbf{I. Single soliton initial data.}
As a first example, we consider the long-time dynamics of one-soliton solutions. In particular, we simulate numerically the time evolution of the gKdV equation \eqref{gkdv} with $F(U) = U + \delta U^3$ and $\delta=0.02$, supplemented with the KdV one-soliton initial datum $u_{\text s}(x)$ in  \eqref{e:kdvsoliton} with $c=0.3$ and $x_0=0$ for which $\no{u_{\text s}}_{H^2}\approx 2.05=\epsilon$. 
While the proximity upper bounds of Theorem \ref{prox-t} exhibit a large growth rate in terms of $\epsilon>1$, a further inspection of the derivation of the relevant estimates (e.g. see inequalities \eqref{IN1} and \eqref{IN2}) suggests that the overall linear growth rate should be compensated by the smallness of the parameter $\delta$ in the example under consideration. The numerical results, illustrated in Figure~\ref{f:2} over a time interval $t\in [0, T]$ with $T\approx 500$, confirm this compensation. Panels (a) and (b) in Figure~\ref{f:2} illustrate snapshots at specific times comparing the profiles of the solutions. They show strong agreement between the solution dynamics of the integrable KdV equation and its non-integrable perturbation, even at large times.  Panel (c)  demonstrates the moderate linear growth of the distances $\no{\Delta(t)}_{\mathcal{X}}$ for $\mathcal{X}=H_x^1$, $L_x^2$ and $L_x^{\infty}$.

Similar numerical results are obtained for the second example, where we simulate the time evolution of the gKdV equation \eqref{gkdv} with $F(U) = U^2 + \delta U^4$ and $\delta=0.01$. This evolution is initialized via the mKdV one-soliton solution $U_{\text s}(x)$ from \eqref{e:kdvsoliton} with $c=0.2$ and $x_0=0$. For this initial datum, $\no{U_{\text s}}_{H^2}\approx 2.41=\epsilon$. The simulation is performed over a time interval $t\in [0, T]$ with $T\approx 500$. As shown in panels (d) and (e) of Figure~\ref{f:2}, we once again observe strong agreement between the solution dynamics of the integrable mKdV equation and its non-integrable perturbation, even at large times, together with the moderate linear growth of the distances $\no{\Delta(t)}_{\mathcal{X}}$ for $\mathcal{X}=H_x^1$, $L_x^2$ and $L_x^{\infty}$,  shown in panel (f).

\begin{figure}[b!]
\centering
\begin{tabular}{cc}
\textbf{(a)}\hspace{5.0cm}\textbf{(b)}\\	
\includegraphics[width=6cm]{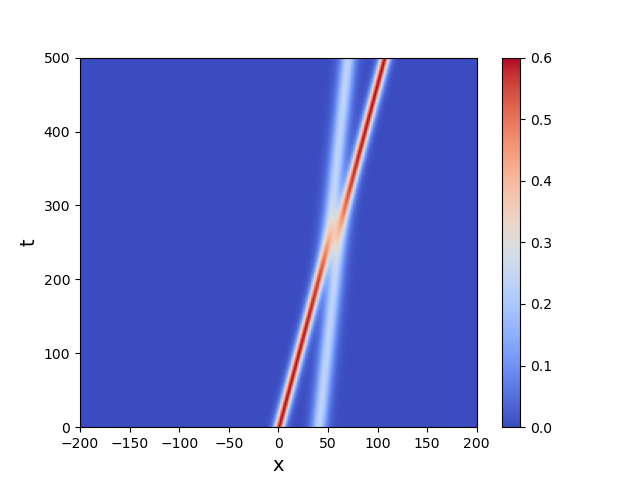}
\includegraphics[width=6cm]{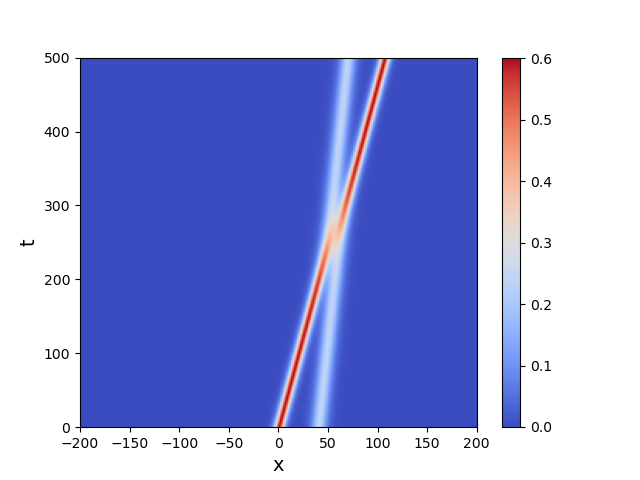}
\\
\textbf{(c)}\hspace{5.0cm}\textbf{(d)}\\	
\includegraphics[width=6cm]{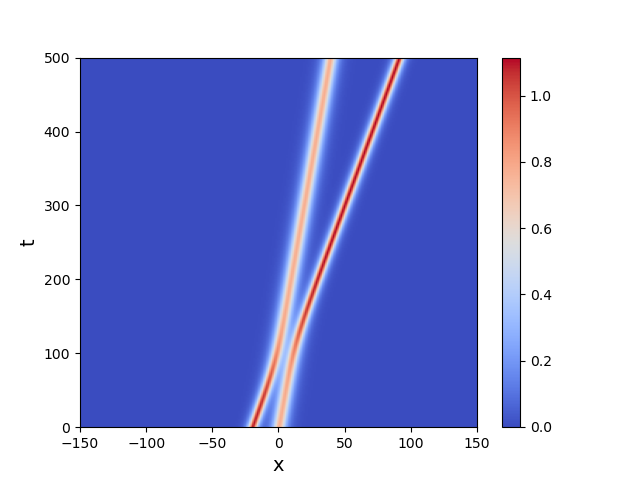}
\includegraphics[width=6cm]{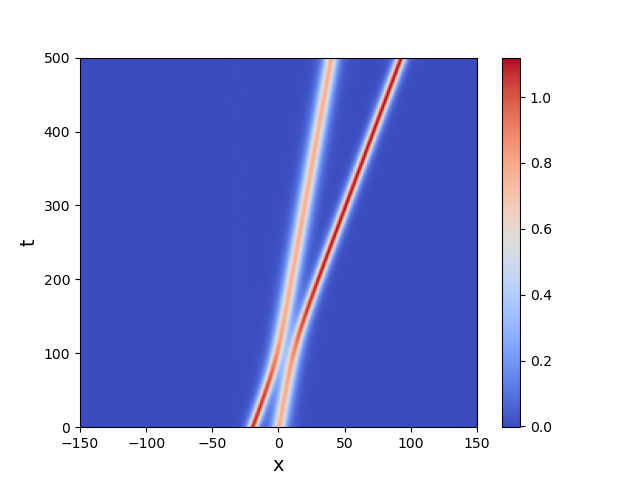}
\end{tabular}
\caption{(a) Contour plot of the (analytical) solution of the KdV equation, supplemented with the KdV two-soliton initial datum \eqref{2sol} with $c_1=0.08$, $c_2=0.2$, $x_{01}=40$ and $x_{02}=0$. (b)  Contour plot of the solution of the gKdV equation \eqref{gkdv} with $F(U)=U+\delta U^3$ and $\delta=0.02$, supplemented with the KdV two-soliton initial datum \eqref{2sol} with $c_1=0.08$, $c_2=0.2$, $x_{01}=40$ and $x_{02}=0$. (c) Contour plot of the (analytical) solution of the mKdV equation, supplemented with the mKdV two-soliton initial datum \eqref{2sol2} with $c_1=0.2$, $c_2=0.1$, $x_{01}=-20$ and $x_{02}=0$. (d) Contour plot of the solution of the  gKdV equation \eqref{gkdv} with $F(U)=U^2+\delta U^4$ and $\delta=0.01$, supplemented with the mKdV two-soliton initial datum \eqref{2sol2} with $c_1=0.2$, $c_2=0.1$, $x_{01}=-20$ and $x_{02}=0$.
}
\label{f:3}
\end{figure}

\begin{figure}[t!]
\centering
\begin{tabular}{cc}
\textbf{(a)}\hspace{5.0cm}\textbf{(b)}\hspace{5.0cm}\textbf{(c)}\\	
\includegraphics[width=5cm]{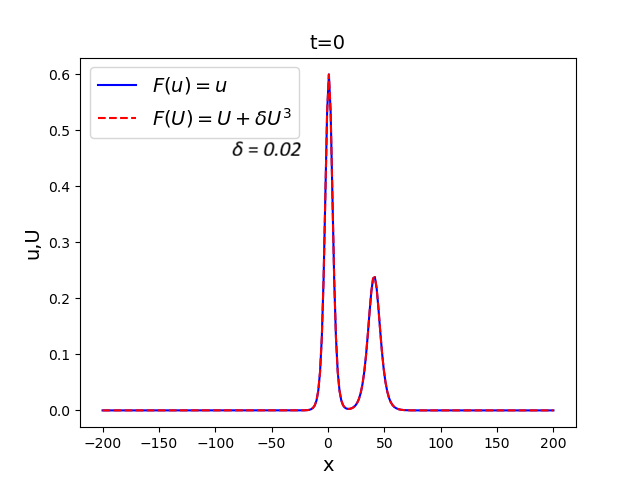}
\includegraphics[width=5cm]{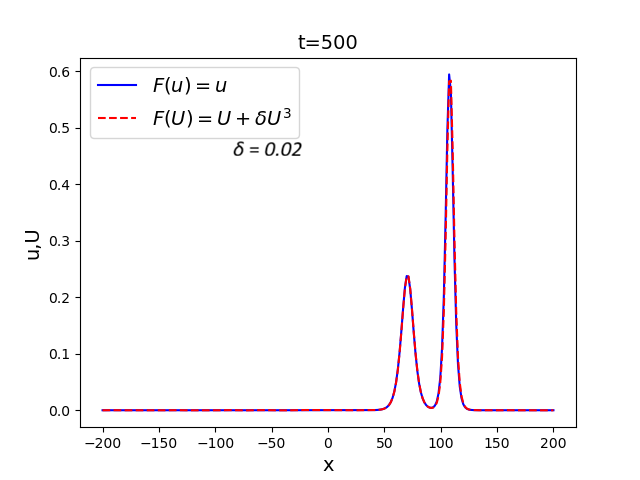}
\includegraphics[width=5cm]{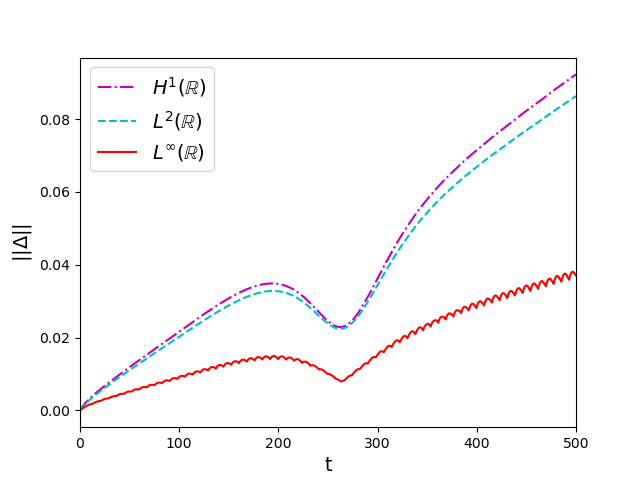}
\\
\textbf{(d)}\hspace{5.0cm}\textbf{(e)}\hspace{5.0cm}\textbf{(f)}\\	
\includegraphics[width=5cm]{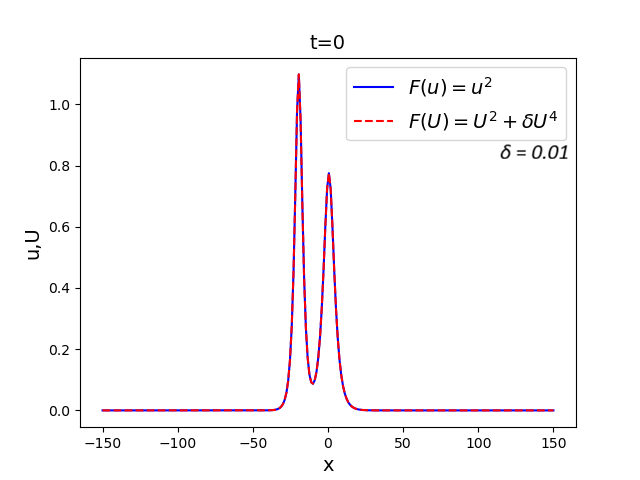}
\includegraphics[width=5cm]{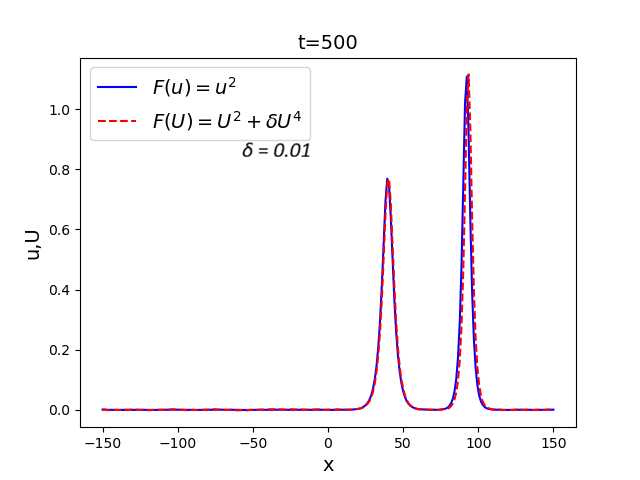}
\includegraphics[width=5cm]{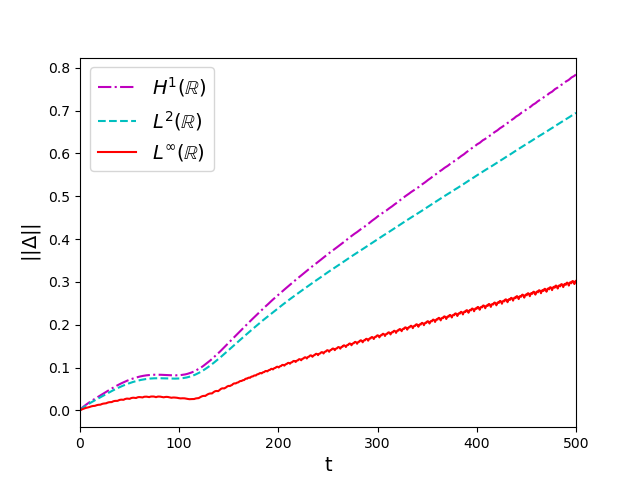}
\end{tabular}
\caption{(a) and (b) Snapshots comparing the profile of the solution of the gKdV equation \eqref{gkdv} with $F(U)=U+\delta U^3$ and $\delta=0.02$, supplemented with the KdV two-soliton initial datum \eqref{2sol} with $c_1=0.08$, $c_2=0.2$, $x_{01}=40$ and $x_{02}=0$, against the profile of the corresponding KdV analytical two-soliton solution. Panel (a) shows the initial profiles and panel (b) the solutions at $t=500$. (c) The corresponding evolution of the norms $\no{\Delta(t)}_{\mathcal{X}}$ with $\mathcal{X}=H_x^1$, $L_x^2$ and $L_x^{\infty}$. (d) and (e) Snapshots comparing the profile of the solution of the gKdV equation \eqref{gkdv} with $F(U)=U^2+\delta U^4$ and $\delta=0.01$, supplemented with the mKdV two-soliton initial datum \eqref{2sol2} with $c_1=0.2$, $c_2=0.1$, $x_{01}=-20$ and $x_{02}=0$, against the profile of the corresponding mKdV analytical two-soliton solution.  Panel (d) shows the initial profiles and panel (e) the solutions at $t=500$. (f) The corresponding time evolution of the norms $\no{\Delta(t)}_{\mathcal{X}}$ with $\mathcal{X}=H_x^1$, $L_x^2$ and $L_x^{\infty}$.
}
\label{f:4}
\end{figure}

\vspace{2mm}
\noindent
\textbf{II. Two-soliton initial data and soliton collisions.}
Next, we consider the interesting case of soliton collisions and, more specifically, the non-integrable dynamics emerging from initial data corresponding to two-soliton solutions of the integrable KdV and mKdV equations. We first study the time evolution of the gKdV equation \eqref{gkdv} with $F(U) = U + \delta U^3$ and $\delta=0.02$, supplemented with KdV two-soliton initial data given by 
\eee{\label{2sol}
u_{\text{2s}}(x) = \sum_{j=1}^2 3c_j\sech^2\left(\frac{\sqrt{c_j}}{2}(x-x_{0j})\right),
}
where $c_1=0.08$, $c_2=0.2$, $x_{01}=40$ and $x_{02}=0$.  %
The norm of this initial datum is $\no{u_{\text 2s}}_{H^2}\approx 0.84=\epsilon$. Since $\epsilon<1$, according to Theorem \ref{prox-t}, we expect proximal dynamics to persist for long times. We perform the simulation over a time interval $t \in [0, T]$ with $T\simeq 500$. In Figure~\ref{f:3}, panel (a) depicts the contour plot of the analytical two-soliton solution of the KdV equation while panel (b)  depicts the contour plot of the solution of the gKdV equation initiated by the above initial datum. For the same example, in Figure~\ref{f:4} panel (a)  shows the initial profiles and panel (b)  shows the comparison of the profiles at $t=500$.   In line with the theoretical analysis, we observe the  remarkable agreement between the soliton collision dynamics of the integrable KdV equation and its non-integrable perturbation. Regarding the evolution of the distances $\no{\Delta(t)}_{\mathcal{X}}$ shown in panel (c) of Figure~\ref{f:4}, it is interesting to note, once again, the appearance of maxima and minima in their time evolution, like in the case of large data for the power nonlinearity (see part II in Section \ref{41-ss}).  These extrema are due to the slight curvature introduced by the weak nonlinearity of the gKdV equation, such that the distance is minimized precisely at the instant of the soliton collision where the caustics of the solitons cross one another.

The second example is devoted to the time evolution of the weakly perturbed mKdV equation, i.e. the gKdV equation \eqref{gkdv} with $F(U) = U^2 + \delta U^4$ and $\delta=0.01$, supplemented with the mKdV two-soliton initial data given~by 
\eee{\label{2sol2}
U_{\text{2s}}(x) = \sum_{j=1}^2 \sqrt{6c_j}\sech\left(\sqrt{c_j}(x-x_{0j})\right),
}
where $c_1=0.1$, $c_2=0.2$, $x_{01}=-20$ and $x_{02}=0$. The norm of this initial datum is $\no{U_{\text 2s}}_{H^2}\approx 3.13=\epsilon$. Since $\epsilon>1$, according to Theorem \ref{prox-t}, we expect the dynamics between the compared equations to diverge at a faster rate which should, however, be compensated again by the smallness of the strength of the weak perturbation $\delta$.

The results of the simulation for $t\in [0, T]$ with $T\approx 500$ are depicted in panels (c) and (d) of Figure~\ref{f:3} and panels (d), (e) and (f) of Figure~\ref{f:4}. In Figure~\ref{f:3},  panel (c) depicts the contour plot of the corresponding analytical two-soliton solution of the mKdV equation while panel (d) depicts the contour plot of the  solution of the gKdV equation initiated by the above initial datum. For the same example, in Figure~\ref{f:4} panel (d) shows the initial profiles and panel (e)  shows the comparison of the profiles at $t=500$.

Indeed, panel (f) of Figure~\ref{f:4}, which portrays the evolution of the norms $\no{\Delta(t)}_{\mathcal{X}}$, shows the divergence of norms at a significantly faster rate than in the case of small initial data of the first example for the weakly perturbed KdV equation. Specifically, comparing with panel (c) of Figure~\ref{f:4}, we observe that the growth rate for the large data is of $\mathcal{O}(10^{-1})$. This rate, although notably larger than the rate $\mathcal{O}(10^{-2})$  observed in the case of small initial data, is still moderate due to the small strength of the perturbation $\delta$ of the mKdV equation. This moderate growth still guarantees the close proximity of the collision dynamics between the KdV equation and the perturbed mKdV equation of the example, as shown in panels (d) and (e) of Figure~\ref{f:3} and the snapshots of the profiles of the solutions shown in panels (d) and (e) of Figure~\ref{f:4}. Returning to  panel (f) of Figure~\ref{f:4} for the distances $\no{\Delta(t)}_{\mathcal{X}}$, we again notice the appearance of the extrema, at shorter times. This observation is in accordance with the fact that the weak perturbation $\delta U^4$ introduces a larger curvature effect in the caustics of the solitons, similar to the case of $k=4$ for the power nonlinearity (see part II of Section \ref{41-ss}). The minimum appears again precisely at the instant of the soliton collisions, as explained in the first example for the pertubed KdV equation.

\section{Conclusion}

An important question in the study of nonlinear dispersive equations concerns the extent to which the dynamics generated by integrable equations persists under non-integrable perturbations and generalizations~\cite{dz2002,t2006,tao2009,deift2017}.
Following a recent investigation in the context of NLS-type equations~\cite{hkmms2024}, in the present work we address the above question within the framework of KdV-type equations. In particular, we consider an important class of generalized, non-integrable KdV (gKdV) equations featuring polynomial nonlinearities, namely the family of equations \eqref{gkdv}, including both standard power-type nonlinearities and weak perturbations of the integrable KdV equation as special cases.

An important first step in our analysis is the derivation of a size estimate for local solutions of the gKdV family~\eqref{gkdv} in the class of Sobolev spaces $H^s(\mathbb{R})$ for any $s\in\mathbb{N}$ with $s>1$. This estimate, which is stated in Theorem \ref{se-t}, is linear in the $H^s(\R)$ norm of the initial data and serves as an analogue of the bound established in~\cite{et2016book} for the  KdV equation. 

The size estimate enables us to proceed to the next step, namely the proof of estimates in the $H^s(\R)$ norm for the distance between solutions of the integrable KdV equation and its non-integrable gKdV counterparts.
These proximity estimates, which are stated in Theorem \ref{prox-t}, predict at most linear in time growth of the deviation between solutions in the $H^s(\R)$ norm, indicating that integrable dynamics originating from small initial data may persist in the non-integrable setting over significant time intervals. 

The numerical simulations presented in Section \ref{s:numerical} nicely illustrate this theoretical prediction. Specifically, for both one-soliton and two-soliton initial data, the observed deviation between the solutions to the KdV and gKdV equations remains small and exhibits linear growth in time, even for durations exceeding the theoretically predicted timescales. This persistence is also observed in the case of weak non-integrable perturbations and, in accordance with the theory, occurs even for initial data of amplitude close to one.

For larger initial data, our theoretical estimates predict deviation at earlier times, and this is indeed illustrated by the numerical simulations. However, in the examples with power-type nonlinearities, the simulations reveal that this deviation arises primarily due to the fact that the KdV soliton, when evolved under the gKdV dynamics, propagates slower than its integrable counterpart, despite exhibiting a robust evolution and preserving its profile up to small radiations. By analyzing the effect of a rotation of the linear caustic of the KdV soliton, we observe that the non-integrable soliton remains remarkably close to the appropriately rescaled integrable KdV soliton. The required scaling parameter depends on the nonlinearity exponent and determines the significantly large time intervals over which proximal dynamics persist.

The theoretical and numerical results of the present work pave the way for interesting potential extensions. On the theoretical side, an important direction is the extension of these results to functional space settings appropriate for studying the persistence of spatially periodic solutions, quasi-periodic solutions, or the potential proximal dynamics originating from non-smooth data. Further theoretical directions concern the case of more general nonlinearities, including for example fractional power nonlinearities. 
Finally, by exploiting complete integrability of the governing model, mathematically precise descriptions of phenomena such as small-dispersion limits, dispersive shock-waves, and soliton gases (statistical ensembles of interacting solitons) have been obtained. 
An interesting future direction is the combination of these descriptions with the proximity results of the present work in order to obtain analogous descriptions for non-integrable models.

\bibliographystyle{myamsalpha}
\bibliography{references.bib}

\newcommand{\etalchar}[1]{$^{#1}$}
\providecommand{\bysame}{\leavevmode\hbox to3em{\hrulefill}\thinspace}
\providecommand{\MR}{\relax\ifhmode\unskip\space\fi MR }
\providecommand{\MRhref}[2]{%
  \href{http://www.ams.org/mathscinet-getitem?mr=#1}{#2}
}
\providecommand{\href}[2]{#2}
\begin{thebibliography}{GGKM67}

\bibitem[AKNS74]{akns1974}
M.~J. Ablowitz, D.~J. Kaup, A.~C. Newell, and H. Segur, \emph{The inverse
  scattering transform-{F}ourier analysis for nonlinear problems}, Studies in
  Appl. Math. \textbf{53} (1974), no.~4, 249--315. \MR{450815}

\bibitem[AS81]{ablowitz1981solitons}
M.~J. Ablowitz and H. Segur, \emph{Solitons and the inverse scattering
  transform}, SIAM, 1981.

\bibitem[Bou93]{b1993-kdv}
J. Bourgain, \emph{Fourier transform restriction phenomena for certain lattice
  subsets and applications to nonlinear evolution equations. {II}. {T}he
  {K}d{V}-equation}, Geom. Funct. Anal. \textbf{3} (1993), no.~3, 209--262.
  \MR{1215780}

\bibitem[BS75]{bs1975}
J.~L. Bona and R. Smith, \emph{The initial-value problem for the {K}orteweg-de
  {V}ries equation}, Philos. Trans. Roy. Soc. London Ser. A \textbf{278}
  (1975), no.~1287, 555--601. \MR{385355}

\bibitem[CKS{\etalchar{+}}03]{ckstt2003}
J. Colliander, M. Keel, G. Staffilani, H. Takaoka, and T. Tao, \emph{Sharp
  global well-posedness for {K}d{V} and modified {K}d{V} on {$\Bbb R$} and
  {$\Bbb T$}}, J. Amer. Math. Soc. \textbf{16} (2003), no.~3, 705--749.
  \MR{1969209}

\bibitem[Dei17]{deift2017}
P. Deift, \emph{Some open problems in random matrix theory and the theory of
  integrable systems. {II}}, SIGMA Symmetry Integrability Geom. Methods Appl.
  \textbf{13} (2017), Paper No. 016, 23. \MR{3622647}

\bibitem[DZ93]{dz1993}
P. Deift and X. Zhou, \emph{A steepest descent method for oscillatory
  {R}iemann-{H}ilbert problems. {A}symptotics for the m{K}d{V} equation}, Ann.
  of Math. (2) \textbf{137} (1993), no.~2, 295--368. \MR{1207209}

\bibitem[DZ02]{dz2002}
P. Deift and X. Zhou, \emph{Perturbation theory for infinite-dimensional
  integrable systems on the line. {A} case study}, Acta Math. \textbf{188}
  (2002), no.~2, 163--262. \MR{1947893}

\bibitem[ET16]{et2016book}
M.~B. Erdogan and N. Tzirakis, \emph{Dispersive {P}artial {D}ifferential
  {E}quations}, London Mathematical Society Student Texts, vol.~86, Cambridge
  University Press, Cambridge, 2016. \MR{3559154}

\bibitem[FFM80]{flaschka1980multiphase}
H. Flaschka, M.~G. Forest, and D. McLaughlin, \emph{Multiphase averaging and
  the inverse spectral solution of the {K}orteweg-de {V}ries equation}, Comm.
  Pure Appl. Math. \textbf{33} (1980), no.~6, 739--784.

\bibitem[GG83]{gear1983second}
J.~A. Gear and R. Grimshaw, \emph{A second-order theory for solitary waves in
  shallow fluids}, The Physics of Fluids \textbf{26} (1983), no.~1, 14--29.

\bibitem[GGKM67]{ggkm1967}
C. Gardner, C. Greene, M. Kruskal, and R. Miura, \emph{Method for solving the
  {K}orteweg-de {V}ries equation}, Phys. Rev. Lett. \textbf{19} (1967), 1095.

\bibitem[GT09]{grunert2009long}
K. Grunert and G. Teschl, \emph{Long-time asymptotics for the {K}orteweg-de
  {V}ries equation via nonlinear steepest descent}, Mathematical Physics,
  Analysis and Geometry \textbf{12} (2009), no.~3, 287--324.

\bibitem[HKM{\etalchar{+}}24]{hkmms2024}
D. Hennig, N.~I. Karachalios, D. Mantzavinos, J. Cuevas-Maraver, and I.~G.
  Stratis, \emph{On the proximity between the wave dynamics of the integrable
  focusing nonlinear {S}chr\"{o}dinger equation and its non-integrable
  generalizations}, J. Differential Equations \textbf{397} (2024), 106--165.
  \MR{4720515}

\bibitem[Kat83]{k1983}
T. Kato, \emph{On the {C}auchy problem for the (generalized) {K}orteweg-de
  {V}ries equation}, Studies in applied mathematics, Adv. Math. Suppl. Stud.,
  vol.~8, Academic Press, New York, 1983, pp.~93--128. \MR{759907}

\bibitem[KPV89]{kpv1989}
C.~E. Kenig, G. Ponce, and L. Vega, \emph{On the (generalized) {K}orteweg-de
  {V}ries equation}, Duke Math. J. \textbf{59} (1989), no.~3, 585--610.
  \MR{1046740}

\bibitem[KPV91]{kpv1991-kdv}
\bysame, \emph{Well-posedness of the initial value problem for the
  {K}orteweg-de {V}ries equation}, J. Amer. Math. Soc. \textbf{4} (1991),
  no.~2, 323--347. \MR{1086966}

\bibitem[KPV96]{kpv1996}
\bysame, \emph{A bilinear estimate with applications to the {K}d{V} equation},
  J. Amer. Math. Soc. \textbf{9} (1996), no.~2, 573--603. \MR{1329387}

\bibitem[Lax68]{l1968}
P.~D. Lax, \emph{Integrals of nonlinear equations of evolution and solitary
  waves}, Comm. Pure Appl. Math. \textbf{21} (1968), 467--490. \MR{235310}

\bibitem[Lax75]{lax1975periodic}
P.~D. Lax, \emph{Periodic solutions of the {K}d{V} equation}, Comm. Pure Appl.
  Math. \textbf{28} (1975), no.~1, 141--188.

\bibitem[Lax76]{lax1976almost}
\bysame, \emph{Almost periodic solutions of the {K}d{V} equation}, SIAM review
  \textbf{18} (1976), no.~3, 351--375.

\bibitem[LP09]{lp2009}
F. Linares and G. Ponce, \emph{Introduction to {N}onlinear {D}ispersive
  {E}quations}, Universitext, Springer, New York, 2009. \MR{2492151}

\bibitem[McK77]{mckean1977stability}
H. McKean, \emph{Stability for the {K}orteweg-de {V}ries equation}, Comm. Pure
  Appl. Math. \textbf{30} (1977), no.~3, 347--353.

\bibitem[McL81]{mclaughlin1981modulations}
D. McLaughlin, \emph{Modulations of {K}d{V} wavetrains}, Physica D: Nonlinear
  Phenomena \textbf{3} (1981), no.~1-2, 335--343.

\bibitem[Miu76]{miura1976korteweg}
R.~M. Miura, \emph{The {K}orteweg-de {V}ries equation: a survey of results},
  SIAM review \textbf{18} (1976), no.~3, 412--459.

\bibitem[Sj{\"o}70]{sjo1970}
A. Sj{\"o}berg, \emph{On the {K}orteweg-de {V}ries equation: existence and
  uniqueness}, J. Math. Anal. Appl. \textbf{29} (1970), 569--579. \MR{410135}

\bibitem[ST76]{st1976}
J.~C. Saut and R. Temam, \emph{Remarks on the {K}orteweg-de {V}ries equation},
  Israel J. Math. \textbf{24} (1976), no.~1, 78--87. \MR{454425}

\bibitem[Sta97]{s1997}
G. Staffilani, \emph{On the generalized {K}orteweg-de {V}ries-type equations},
  Differential Integral Equations \textbf{10} (1997), no.~4, 777--796.
  \MR{1741772}

\bibitem[Tao06]{t2006}
T. Tao, \emph{Nonlinear {D}ispersive {E}quations}, CBMS Regional Conference
  Series in Mathematics, vol. 106, American Mathematical Society, 2006, Local
  and global analysis. \MR{2233925}

\bibitem[Tao09]{tao2009}
\bysame, \emph{Why are solitons stable?}, Bull. Amer. Math. Soc. (N.S.)
  \textbf{46} (2009), no.~1, 1--33. \MR{2457070}

\bibitem[Tem69]{t1969}
R. Temam, \emph{Sur un probl\`eme non lin\'eaire}, J. Math. Pures Appl. (9)
  \textbf{48} (1969), 159--172. \MR{261183}

\bibitem[ZK65]{zk1965}
N. Zabusky and M.~D. Kruskal, \emph{Interaction of solitons in a collisionless
  plasma and the recurrence of initial states}, Phys. Rev. Lett \textbf{15}
  (1965), 240--243.

\end{thebibliography}

\end{document}